\documentclass[12pt]{article}
\usepackage[top=1in,bottom=1in,left=0.5in,right=0.5in]{geometry}
\usepackage{amsmath,amssymb,amsfonts,bbm, mathrsfs}
\usepackage{graphicx}
\usepackage{natbib}
\usepackage{setspace}
\usepackage{algorithmic}
\usepackage{algorithm}
\usepackage{caption}
\newcommand{\Bta}{\mathrm{Beta}}

\newcommand{\argmax}{\mathrm{argmax}}

\newcommand{\indicator}[1]{\mathbbm{1}_{\left[ {#1} \right] }}
\newcommand{\zap}[1]{}

\usepackage{natbib}
 \bibpunct[, ]{(}{)}{,}{a}{}{,}%

\def\FIGURE#1#2#3{%
\begin{centering}
#1
\caption{#2}
\end{centering}
}

\newtheorem{theorem}{Theorem}
\newtheorem{proposition}{Proposition}
\newtheorem{lemma}{Lemma}
\newtheorem{remark}{Remark}

\begin{document}
\title{Exploration vs. Exploitation\\in the Information Filtering Problem}
\author{Xiaoting Zhao and Peter I. Frazier}
\date{\today}
\maketitle
\abstract{We consider information filtering, in which we face a stream of items too voluminous to process by hand (e.g., scientific articles, blog posts, emails), and must rely on a computer system to automatically filter out irrelevant items.  Such systems face the exploration vs. exploitation tradeoff, in which it may be beneficial to present an item despite a low probability of relevance, just to learn about future items with similar content.  We present a Bayesian sequential decision-making model of this problem, show how it may be solved to optimality using a decomposition to a collection of two-armed bandit problems, and show structural results for the optimal policy.  We show that the resulting method is especially useful when facing the cold start problem, i.e., when filtering items for new users without a long history of past interactions.  We then present an application of this information filtering method to a historical dataset from the arXiv.org repository of scientific articles.

}
\section{Introduction}
 We consider the information filtering problem, in which a human (the {\it user}) is tasked with processing a stream of items (e.g., emails, text documents, intelligence information, bug reports, scientific articles). Some of these items are relevant and should be examined in detail, while the rest are irrelevant and can be ignored.  When the stream is too voluminous to be processed by hand, a computer system can be tasked with automatically pre-processing, or {\it filtering}, these items, forwarding some on to the user, and discarding others.  In creating such an information filtering system, we wish to forward as many relevant items as possible, without forwarding too many irrelevant items.

Information filtering systems typically rely on large amounts of historical data to fit statistical models for predicting item relevance (e.g., \cite{AgChPa2011},  \cite{AgarwalZhangMazumder2011}, \cite{ShivaswamyJoachims12b}). When access to historical information is limited, it is difficult to build an effective information filtering system \citep{hanani2001information}.  This makes information filtering difficult to apply for new users, new types of relevance, or when items' characteristics or users' interests are rapidly evolving. This so-called ``cold start'' problem is prevalent in many information filtering systems \citep{Schein02,RubensRecSysHB2010}. 

When historical data is limited, an information filtering system can also learn about item relevance from user-provided feedback (implicitly, through clicks, or explicitly, through ratings) about previously forwarded items.  Moreover, when faced with a sub-stream of items originally predicted to be irrelevant based on limited historical data, an information filtering system might forward a small number of these to the user for feedback, learning with greater certainty their true relevance. Such {\it exploration} of user preferences is useful because, if this sub-stream is revealed to be relevant, future items from that sub-stream can be forwarded.

However, too much exploration will lead to too many irrelevant forwarded items.  Thus, an information filtering system should also put some weight on {\it exploitation}, i.e., forwarding only those items predicted to be relevant with a high degree of certainty.  This tradeoff between exploration and exploitation, which appears in other problem domains including reinforcement learning \citep{SuttonBartoRL98,AuerNearOptimalRegret2010}, approximate dynamic programming \citep{Powell04learningalgorithms,PowellADP2007}, revenue management \citep{ArCa10,BeZe09,den2013simultaneously}, and inventory control \citep{LaPo99,DiPuBi02},
is also important for understanding the information filtering problem in regimes with little historical data.

In this paper, we propose and analyze a mathematical model of the information filtering problem, formulating it as a stochastic control problem using Bayesian statistics and stochastic dynamic programming. Our analysis provides insight into the exploration vs. exploitation tradeoff in information filtering, and more generally into the structure of the optimal information filtering strategy. In comparison with a myopic ``pure exploitation'' strategy, we show that the optimal filtering strategy forwards every item that the myopic strategy forwards, and potentially forwards additional items that the myopic strategy would not forward.  Moreover, this willingness to forward additional items, i.e., to explore, is largest when the number of items on which we have relevance feedback is small, and decreases as this number of items with feedback grows larger.  In the limit as the number of items with feedback grows to infinity, this willingness to forward additional items vanishes, and the decisions of the optimal strategy match those of the myopic strategy.

We additionally provide an efficient method for computing the optimal information filtering strategy. While the optimal strategy is the solution to a partially observable Markov Decision Processes, and thus can be computed, at least conceptually, using dynamic programming \citep{FrWeorDp}, the curse of dimensionality prevents directly computing this solution in practice. To circumvent this issue, we show that the problem can be decomposed into a collection of much smaller dynamic programs, which can be solved efficiently. Indeed, each smaller dynamic program is a two-armed bandit problem, with one unknown Bernoulli arm and one known arm, and can be solved directly, or using methods from the extensive literature on two-armed bandits \citep{Bellman1956, GiJo74, Gittins79, Wh80, BeFr85, KaVe87, GiGlWe11}. 

Our model can be seen as a Bayesian contextual bandit problem with two arms (forward and discard), and a particular structure (discussed in Section~\ref{sec:model}) for the relationship between context and reward. While finding an optimal policy for general Bayesian contextual bandits is challenging \citep{LaZh07, MayKordaLeeLeslie2012, AgGo2012}, the special structure that we assume allows us to compute an optimal policy tractably.  

This work is motivated by an information filtering system we are building for the online repository of scientific articles, arXiv.org \citep{arxiv}.  By December 2014, arXiv.org had accumulated over 1 million full-text articles, was receiving an additional 7000+ new articles each month, and was distributing about 1 million downloads weekly to 400,000 unique users \citep{Ginsparg2011, VanNoorden2014}. This massive stream of articles creates a challenge for researchers who wish to keep abreast of those new articles relevant to their research.  In Section~\ref{sec:simulation}, we present results from applying our proposed algorithm to arXiv.org, and show that exploring in an optimal way can provide significant value over three benchmark policies: pure exploitation, Thompson sampling \citep{Thompson1933, AgGo2011, ChapelleLi2011, RussoVanRoy2014}, and upper confidence bound (UCB) \citep{LaRo85, AuCeFi02, Kaufmann2012}.



There are many works in the active learning and information retrieval communities researching information filtering systems. For an overview, see \cite{RubensRecSysHB2010}, \cite{IR_manning2008}, and \cite{AdTu05}. The earliest stage of research focuses on cost/credit of delivering relevant/irrelevant documents by various classifiers, including support vector machines \citep{Joachims98textcategorization},  inference networks \citep{Callan96documentfiltering}, and maximizing historical data likelihood \citep{Lafferty01documentlanguage,ZhCa01}. In almost all of this work, the future benefit of reducing uncertainty through exploration is ignored, and a pure exploitation policy is used.

The most related work to our approach is \cite{ZhXuCa03}, which studies the exploration vs. exploitation tradeoff in information retrieval using a Bayesian decision-theoretic model.
Using Bayesian logistic regression to measure model quality, that previous work quantifies the one-step value of information associated with observing feedback on an item, and computes this approximately using Monte Carlo.
While our approach shares conceptual similarities, we consider a different model, in which we observe an item's category rather than its score.
Additionally, our work goes beyond this previous work by relaxing the one-step assumption, and providing a policy that is optimal for multiple steps into the future.
Thus, our work can be seen as analogous to similar efforts to go beyond one-step optimality in related Bayesian sequential decision-making problems \citep{XieFrazier2011,JeFrSz11}.

Among other related work, \cite{AgarwalChenElango2009}, \cite{RaKlJo08}, \cite{YuBrKlJo09} and \cite{Hofmann2013} study multi-armed bandit methods in recommender systems,
and \cite{Xu:2008:NPR:1390334.1390408} considers the problem of choosing results to an initial user query, so as to best improve later search results.
While both lines of research study the exploration vs. exploitation tradeoff in information retrieval, neither directly considers the information filtering problem. \cite{Shani2005} introduces the concept of Markov decision processes when modeling recommender systems while \cite{ LeRuMa2013} applies a sequential event prediction technique to recommender systems. Both works focus on the fact that recommendation is a sequential decision problem, where the revenue earned in the current period may depend on more than just the most recent action.

Below, in Section~\ref{sec:model}, we formulate the information filtering problem as a stochastic control problem. In Section~\ref{sec:solution}, we provide an efficient solution by decomposing the original problem into multiple sub-problems that can be solved efficiently. We then show structural results: the optimal policy always forwards at least those items forwarded by a pure exploitation policy, and is a threshold policy whose threshold is non-decreasing in the total number of observed items. We also relate these structural results to known properties of two-armed bandit problems. Lastly, we present experimental results in Section~\ref{sec:simulation}, using both idealized Monte Carlo simulations and trace-driven simulations with historical data from arXiv.org.


\section{Mathematical Model}
\label{sec:model}
Each item arriving from our information stream is labeled with (exactly) one of $k$ categories, and we let $X_n \in \{1, \ldots, k\}$ be the category of the $n$th item in our steam. This category is observable by our information filtering algorithm.  In the application to arXiv.org that we present in Section~\ref{sec:simulation}, we describe a setting where the category is provided by a human (the author) who submits the item to the stream. This category could also be obtained automatically by a machine learning algorithm from item contents.  We model the sequence of random variables $(X_n: n=1,2,\ldots)$ as being independent and identically distributed, and we let $p_x = P(X_n = x)>0$ for $x=1,\ldots,k$. 

In our model and analysis, we focus on a single user, and then in implementation we apply the resulting algorithm separately for each user. Fixing this user, we model each category as having associated with it some latent unobservable value $\theta_x\in[0,1]$, which is the probability that the user under consideration would find an item from this category to be relevant, if it were forwarded to her/him.  We let $\theta = [\theta_1,\ldots,\theta_k]$.  Our model assumes that $\theta_x$ remains static over time. This focus on a single user is in contrast with much of the work on collaborative filtering, and is motivated by the design requirements of the information filtering system we are building for arXiv.org, where concerns about fairness to authors make it especially important to avoid cascades and the Matthew effect \citep{EasleyKleinberg2010}, in which popular items become more popular, irrespective of quality.

We model each $\theta_x$ as having been drawn independently for each $x$ from a Bayesian prior probability distribution, which is beta-distributed.  We let $\alpha_{0x},\beta_{0x}$ be the two parameters of this distribution, so that $\theta_x \sim \mathrm{Beta}(\alpha_{0x},\beta_{0x})$.  In Section~\ref{sec:simulation} we provide a method for estimating the parameters of this prior probability distribution from historical data. 

For each item in the stream, our information filtering algorithm then decides whether to forward this item to the user, or to discard it.
We let $U_{n} \in \{0,1\}$ represent the decision made for the $n$th item, where $U_{n}$ is $1$ if the algorithm forwards this item, and $0$ if it discards it.
Each forwarded item is then seen by the user, who provides feedback on its relevance in the form of a Bernoulli random variable $Y_{n}$.  In the application presented in Section~\ref{sec:simulation}, the feedback $Y_{n}$ is provided implicitly, by whether or not the user clicked to view the forwarded item.  In other applications, this feedback might be provided via an explicit rating inputted by the user.  The sequence of random variables $(Y_{n} : n=1,2,\ldots)$ are conditionally independent given $\theta_x$, with $P(Y_{n}=1 | X_{1:n}, U_{1:n}, \theta) = \theta_{X_n}$.
No feedback is provided on discarded items, and so $Y_n$ is observed if and only if $U_n=1$.

The decision of whether or not to forward the item may be made based only on the information available from previous forwarding decisions.  That is, we require $U_{n}$ to depend only on the current item's category $X_n$, and the history $H_{n-1} = (X_{\ell},U_{\ell},U_{\ell}Y_{\ell} : \ell \le n-1)$.  In this definition of the history, we emphasize that $Y_{\ell}$ is only observable if $U_{\ell}=1$.

In our model, we pay an explicit cost $c$ for each item forwarded to the user, which models the cost of the user's time, and a reward of $1$ for each relevant item forwarded.
Thus, the total reward resulting from forwarding an item is $U_n (Y_n - c)$.

At each time step $n$, there is a probability $\gamma$ that the user will remain engaged through the next time step $n+1$, and a probability $1-\gamma$ that the user will abandon the system and never return, so $N$ follows a geometric distribution with parameter $1-\gamma$. We have modeled the items as arriving in discrete time. Our approach can also be easily adapted to a continuous-time setting, where documents arrive according to a Poisson process and the user has an exponential lifetime in the system. We discuss this further in Section~\ref{sec:simulation} in the context of arXiv.org along with estimation of $\gamma$ and validation of our modelling assumptions. 

\begin{figure}
\FIGURE
  {\includegraphics*[width=4.5in, height=2.5in]{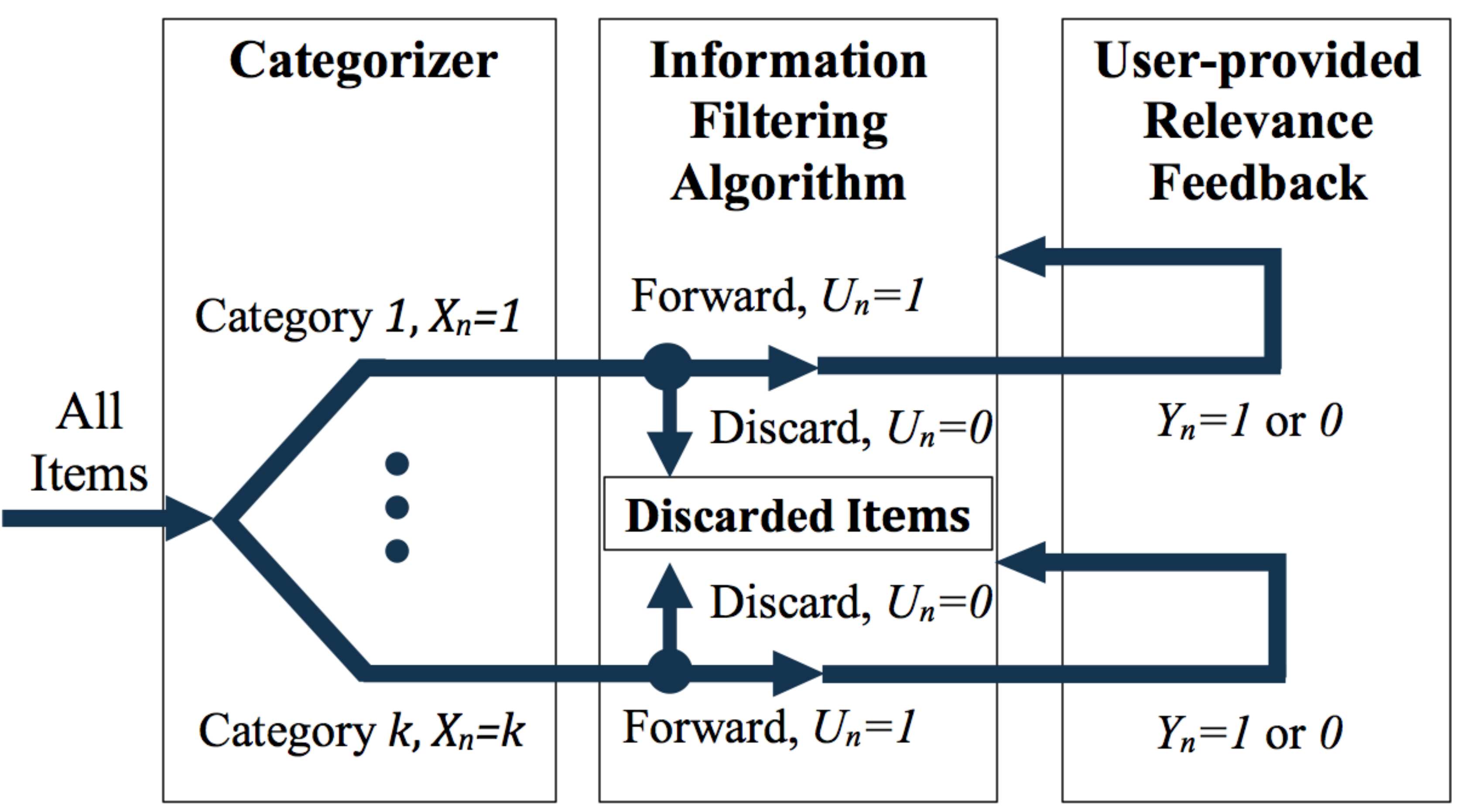}}
  {Schematic of the information filtering problem.  Arriving items are categorized into one of $k$ categories, and then are forwarded or discarded by an information filtering algorithm.  This algorithm uses feedback on forwarded items to improve later forwarding decisions.\label{fig:schematic} }{}
\end{figure}

In Section~\ref{sec:simulation}, we discuss estimation of $\gamma$ in our application to arXiv.org, and validation of our assumption. 
Our goal is to design an algorithm for making forwarding decisions that maximizes the expectation of the cumulative reward, $\sum_{n=1}^N U_n(Y_n-c)$.
Our model is summarized by Figure~\ref{fig:schematic}.

To formalize this as a stochastic control problem, we define a policy $\pi$ as a sequence of functions, $\pi=(\pi_1,\pi_2,\ldots)$, where each $\pi_n : (\{1,\ldots,k\} \times \{0,1\}\times\{0,1\})^{n-1} \times \{1,\ldots,k\} \mapsto \{0,1\}$ maps histories onto actions.  We let $\Pi$ be the space of all such policies.
For each $\pi\in\Pi$, we define $P^\pi$ to be the measure under which $U_{n+1} = \pi_{n+1}(H_n,X_{n+1})$ almost surely for each $n$, and we let $E^\pi$ represent the expectation taken with respect to this measure.
Our goal is then to solve
\begin{equation}
  \sup_{\pi\in\Pi} E^\pi\left[\sum_{n=1}^N U_n(Y_n-c)\right].
  \label{eq:overall-problem}
\end{equation}

\section{Solution Method and Structural Results}
\label{sec:solution}
While \eqref{eq:overall-problem} is a stochastic control problem, and could be solved via dynamic programming, the size of the state space of this dynamic program grows exponentially in $k$.  This is the so-called ``curse of dimensionality'' \citep{PowellADP2007}.  To circumvent this issue, we decompose the problem into a collection of single-category problems, each of which is a two-armed bandit problem and can be solved efficiently using dynamic programming.  This is accomplished in Section~\ref{sec:decomposition}, which performs the decomposition, and Section~\ref{sec:undiscounted-to-discounted}, which converts each single-category problem from its original form (undiscounted random horizon) to an easier-to-solve form (discounted infinite-horizon). We then provide the dynamic programming equations for these single category problems in Section~\ref{sec:single-category}, which we then use to show structural results in Section~\ref{sec:structural}. Sections~\ref{sec:undiscounted-to-discounted} and~\ref{sec:single-category} follow standard arguments, but are included to support results in Section~\ref{sec:structural}. 

\subsection{Decomposition into single-category subproblems}
\label{sec:decomposition}
To decompose the problem \eqref{eq:overall-problem} into a number of easily solved single-category problems, we first introduce some additional notation.
First, we let $n_{\ell x} = \inf\left\{ n: \sum_{i=1}^n \indicator{X_i = x} = \ell \right\}$ be the index, in the overall stream of items, of the $\ell$th item from category $x$.  
We then define $U_{\ell x} = U_{n_{\ell x}}$ and $Y_{\ell x} = Y_{n_{\ell x}}$ to be the forwarding decision and relevance, respectively, of the $\ell$th item from category $x$.  Finally, we let $N_x = \sup\left\{ \ell: n_{\ell x} \le N\right\}$ be the number of items from category $x$ that arrive before the user leaves the system at step $N$.
We let $H_{nx} = (U_{\ell x},U_{\ell x}Y_{\ell x} : \ell \le n)$ be the history of forwarding decisions and relevance feedback for items from category $x$.

For each category $x$, we then define a class of policies for making forwarding decisions for items {\it from that category only}, based only on the portion of the history arising from items in that category.  Formally, we define a {\it single-category policy} for a given category $x$ to be a sequence of functions, $\pi^{(x)} = (\pi^{(x)}_1,\pi^{(x)}_2,\ldots)$, where each function
$\pi^{(x)}_{n+1}: \{0,1\}^{2n} \mapsto \{0,1\}$ maps $H_{nx}$ onto $U_{n+1,x}$.
We let $\Pi^{(x)}$ be the space of all such policies.  For each $\pi^{(x)}\in\Pi^{(x)}$, we let
$P^{\pi^{(x)}}$ be the probability measure under which $U_{n+1,x}=\pi^{(x)}_{n+1}(H_{nx})$ for each $n$, and we let $E^{\pi^{(x)}}$ be the expectation with respect to this measure.

We will write the value of the overall problem \eqref{eq:overall-problem} in terms of the sum of the values of the solutions to forwarding problems for individual categories,
\begin{equation}
  \label{eq:single-category-problem}
  \sup_{\pi^{(x)}\in\Pi^{(x)}} E^{\pi^{(x)}}\left[\sum_{n=1}^{N_x} U_{nx}(Y_{nx}-c)\right].
\end{equation}
We do this in the following theorem. All proofs may be found in the appendix.

\begin{theorem}
  \label{t:decomposition}
  \begin{equation*}
  \sup_{\pi\in\Pi} E^\pi\left[\sum_{n=1}^N U_n(Y_n-c)\right]
  = \sum_{x=1}^k
  \sup_{\pi^{(x)}\in\Pi^{(x)}} E^{\pi^{(x)}}\left[\sum_{n=1}^{N_x} U_{nx}(Y_{nx}-c)\right].
  \end{equation*}

  Moreover, if $\pi^{(x),*}$ attains the supremum in  \eqref{eq:single-category-problem} for each $x$, and if $\pi^*$ is the policy constructed by setting
  $\pi^*_{n+1}(H_n,X_{n+1}) = \pi^{(X_{n+1}),*}_{n+1}(H_{nx})$
  for each $n$,
  then $\pi^*$ attains the supremum in  \eqref{eq:overall-problem}.
\end{theorem}

This theorem gives us a way to construct the optimal policy for the multi-category problem from the solutions to single-category problems: to make a forwarding decision for a new item in the multi-category problem, we identify the category $X_{n+1}$ of that item, and then perform the forwarding decision that would have been made by the optimal policy for the single-category problem for that category.

We give intuition behind the proof of this theorem here.  Because the objective is additive across categories, and because the prior on categories' relevance are independent of each other, relevance feedback from one category gives no information about the relevance of other categories.
Thus, when considering whether or not to forward an item for a particular category, it is sufficient to consider the history of observations from that category alone.

\subsection{Conversion of finite-horizon single-category subproblems into infinite-horizon two-armed bandits}
\label{sec:undiscounted-to-discounted}
Below, it will be useful to transform this single-category problem~\eqref{eq:single-category-problem} from a finite-horizon undiscounted problem, with a random horizon, into an infinite-horizon discounted problem. This transformation will make clear that the single subproblems are two-armed bandit problems. In performing this transformation, we use a standard geometric killing approach.

We first note that the number of items available for forwarding in a particular category has a geometric distribution. 
\begin{remark}
  \label{r:geometric}
  $N_x \sim \mathrm{Geometric}\left(1-\gamma_x \right)$, where $\gamma_x = \frac{p_x \gamma}{p_x \gamma + 1-\gamma}$.
  Here, by $\mathrm{Geometric}(1-q)$, we mean the probability distribution supported on $\{0,1,2,\ldots\}$ that assigns probability mass $(1-q)q^n$ to integer $n$.
\end{remark}

Using this remark, the following lemma writes the performance of any policy $\pi^{(x)}$ as an infinite-horizon discounted sum.
\begin{lemma}
  \label{l:discount}
  For each policy $\pi^{(x)} \in \Pi^{(x)}$, we have
\begin{equation*}
  E^{\pi^{(x)}}\left[\sum_{n=1}^{N_x} U_{nx}(Y_{nx}-c)\right]
  =\gamma_x E^{\pi^{(x)}}\left[\sum_{n=1}^\infty \gamma_x^{n-1} U_{nx}(Y_{nx}-c)\right],
\end{equation*}
\end{lemma}

This lemma shows that we can find an optimal policy for the single category problem \eqref{eq:single-category-problem} by solving the stochastic control problem,
\begin{equation}
  \label{eq:single-category-discounted-problem}
  \sup_{\pi^{(x)}\in \Pi^{(x)}} E^{\pi^{(x)}}\left[\sum_{n=1}^\infty \gamma_x^{n-1} U_{nx}(Y_{nx}-c)\right].
\end{equation}

Any policy that is optimal for this problem is optimal for \eqref{eq:single-category-problem}.
We have dropped the strictly positive constant $\gamma_x$ when stating this stochastic control problem, since this constant does not affect the ordering of the policies.

Equation~\eqref{eq:single-category-discounted-problem} is a two-armed bandit, with one unknown Bernoulli arm (giving rewards of $1-c$ or $-c$) and one known arm (giving rewards of 0).
By adding an additional reward $c$ in all periods, this problem can be seen to be equivalent to a more conventional two-armed bandit, where the unknown Bernoulli arm gives rewards of $0$ or $1$, and the known arm gives rewards of $c$.
Both two-armed and Bernoulli bandits have been studied extensively \citep{Bellman1956, BeFr85}, and a variety of methods have been proposed for computing optimal policies
\cite{KaVe87,GiGlWe11}.
One may use one of these techniques, or simply solve \eqref{eq:single-category-discounted-problem} directly using dynamic programming, as described in the Appendix.

After Section~\ref{sec:single-category} summarizes the dynamic programming equations for \eqref{eq:single-category-discounted-problem}, Section~\ref{sec:structural} will prove structural results that apply beyond the forwarding problem more generally to Bernoulli bandits.  These results do not, to the best of our knowledge, appear previously in the literature.  We will also relate these new results to related results on Bernoulli bandits from \cite{BeFr85}, and for bandit rewards in location-scale families from \cite{GiGlWe11}.

\subsection{Dynamic programming equations for the single-category subproblem}
\label{sec:single-category}
We now provide the dynamic programming equations for the infinite-horizon discounted version of the single-category problem \eqref{eq:single-category-discounted-problem}, 
which we use in Section~\ref{sec:structural} to derive novel structural results. 
As noted above, the same optimal policy is optimal for the finite-horizon undiscounted problem \eqref{eq:single-category-problem}, and can be used with Theorem~\ref{t:decomposition} to provide an optimal policy for the original problem \eqref{eq:overall-problem}. 

First, the conditional distribution of $\theta_x$ given history $H_{nx}$ is
\begin{equation*}
  \theta_x ~|~H_{nx} ~\sim \text{Beta} \left( \alpha_{nx}, \beta_{nx} \right),
\end{equation*}
where $\alpha_{nx}=\alpha_{0x}+\sum^{n}_{i=1}U_{ix} ~Y_{ix}$
is the sum of $\alpha_{0x}$ and the number of relevant items forwarded,
and
$\beta_{nx}=\beta_{0x}+\sum^n_{i=1}(1-Y_{ix}) U_{ix}$
is the sum of $\beta_{0x}$ and the number of irrelevant items forwarded.
This follows from standard results from Bayesian statistics on conjugate priors for Bernoulli observations.  For details, see, e.g., \cite{Dg04}.  It will also be convenient to introduce the notation $\mu(\alpha,\beta) = \alpha/(\alpha+\beta)$ to refer to the mean of a $\mathrm{Beta}(\alpha,\beta)$ distribution.

Moreover, since $Y_{\ell x}$ are conditionally i.i.d. given $\theta_x$, the conditional distribution of the sequence $(Y_{\ell x} : \ell >n)$ given the history $H_{nx}$ is completely determined by $\alpha_{nx},\beta_{nx}$.
Thus, the solution to \eqref{eq:single-category-discounted-problem} will be given by a dynamic program whose state space includes all possible values of $(\alpha_{nx},\beta_{nx})$.

We briefly review this use of dynamic programming, providing definitions and notations that will be used below when presenting and proving structural results.
For any scalar real numbers $\alpha,\beta>0$, we define the value function
\begin{equation}
\label{eq:valueFunction}
  V_x(\alpha,\beta) = \sup_{\pi^{(x)} \in \Pi^{(x)}} E^{\pi^{(x)}}\left[\sum_{n=1}^\infty \gamma_x^{n-1} U_{nx} (Y_{nx}-c) \mid \theta_x \sim \mathrm{Beta}(\alpha,\beta)\right].
\end{equation}

The value function satisfies Bellman's equation
\begin{equation}
\label{eq:Bellman}
V_x(\alpha, \beta)= \max \{ Q(\alpha, \beta, 0), ~Q(\alpha, \beta, 1)\},
\end{equation}
where we define the Q-factor $Q: (0, \infty)^2 \times \{0,1\} \mapsto \mathbb{R}$ by
\begin{align}
Q_x(\alpha, \beta, 0)&=\gamma_x V_x(\alpha, \beta), \label{eq:Q0} \\
Q_x(\alpha, \beta, 1)&=E[Y_1-c+\gamma_x V_x(\alpha_{1x}, \beta_{1x})|U_{1x}=1, \alpha_{0x}=\alpha, \beta_{0x}=\beta]\\
&=\mu-c+\gamma_x \big[\mu V_x(\alpha+1, \beta)+(1-\mu)V_x(\alpha, \beta+1) \big],\label{eq:Q1}
\end{align}
where $\mu=\mu(\alpha,\beta) = \alpha/(\alpha+\beta)$. Here, $Q_x(\alpha, \beta, 0)$ and $Q(\alpha, \beta, 1)$ are the value of discarding and forwarding the item, and proceeding optimally thereafter, when the posterior on $\theta_x$ is $\mathrm{Beta}(\alpha, \beta)$. 

An optimal policy is then any whose decisions attain the maximum in this recursion, breaking ties arbitrarily.  That is, an optimal policy is any for which
\begin{equation*}
  U_{n+1,x} \in \argmax_{u=0,1} Q_x(\alpha_{nx},\beta_{nx},u).
\end{equation*}

Thus, if we are able to compute the value function, from it we may compute the Q-factors, and then compute an optimal policy.  While the set of possible values for $(\alpha_{nx},\beta_{nx})$, and thus the size of the state space that must be considered, is countably infinite, preventing exact computation of the value function, we describe below in the Appendix a truncation method for computing upper and lower bounds on $V_x$, from which an approximation with explicit error bounds can be computed.  The error in this approximation vanishes as the level of truncation grows.

\subsection{Structural Results}
\label{sec:structural}
In this section, we prove structural results that provide insight into the behavior of the optimal policy (that the optimal policy is a threshold policy, in Theorem~\ref{th:optimalPolicy}, and behavior of this threshold, in Theorem~\ref{th:muStar_structure}), yield computational benefits by reducing the amount of storage needed to implement the optimal policy (Theorem~\ref{th:optimalPolicy}), and that form the foundation for computing approximations with explicit error bounds
(Propositions~\ref{prop:lowerBound} and~\ref{prop:upperBound}). Similar properties are known to hold for location-scale families, and our work shows that they also hold in the Beta-Bernoulli setting, which is more relevant in the ``click-based" forwarding problem studied here.

When solving the single-category sub-problem, we first notice that the Q-factor for the {\it discard} decision, $Q_x(\alpha, \beta, 0)$, has the following structure, stated in Remark~\ref{re:notForward}. 
\begin{remark} 
\label{re:notForward}
It is optimal to discard at $(\alpha, \beta) \in (0, \infty)^2$ if and only if $V_x(\alpha, \beta)=0$. Equivalently, $Q_x(\alpha, \beta, 0)=0$. 
\end{remark}

Remark~\ref{re:notForward} implies that we may rewrite Bellman's equation \eqref{eq:Bellman} as 
$V_x(\alpha, \beta)= \max \{ 0, ~Q(\alpha, \beta, 1)\}$.
Next, we provide a lower and upper bound on the value function, and use them below to prove convergence of the value function in Proposition~\ref{prop:convergence}.  We also use them when describing a truncation method with explicit error bounds that computes the solution to the single-category sub-problem in the Appendix.

Since the value of any policy provides a lower bound on the value function, we consider a policy that ignores feedback, forwards all items if $\mu(\alpha_{0x}, \beta_{0}) \geq c$, and discards all items otherwise. The value of this policy is easy to compute, and so provides a convenient lower bound.  This is the basis for Proposition~\ref{prop:lowerBound}.
\begin{proposition}
\label{prop:lowerBound}
Define $V_x^L(\alpha, \beta)= \frac 1{1-\gamma_x}\max \left \{0, ~\mu(\alpha,\beta)-c \right \}$ for any $\alpha \in (0, \infty)$ and $\beta \in (0, \infty)$. Then, we have the lower bound, $V_x(\alpha,\beta)\geq V_x^L(\alpha, \beta)$.
\end{proposition}

To obtain an upper bound on $V_x(\alpha,\beta)$, we consider an environment in which the true value of $\theta_x$ is revealed. The value of an optimal policy in this environment provides an upper bound on the value function.
This is the basis for Proposition~\ref{prop:upperBound}.

\begin{proposition}
\label{prop:upperBound}
 Define $V_x^U(\alpha, \beta) = \frac1{1-\gamma_x}E[\max\{0, \theta_x-c\}]$ for any $\alpha \in (0, \infty)$ and $\beta \in (0, \infty)$. Then, we have the upper bound $V_x^{U}(\alpha,\beta) \geq V_x(\alpha,\beta).$
\end{proposition}

Combining Proposition~\ref{prop:lowerBound} and Proposition~\ref{prop:upperBound}, we show convergence of the value function as the prior converges to one in which we are certain about $\theta_x$, as stated below in Proposition~\ref{prop:convergence}. This proposition is used in the proofs of Lemma~\ref{l:V_nonDec} and~\ref{l:V_convex}, and of Theorem~\ref{th:muStar_structure}.

\begin{proposition}
\label{prop:convergence}
Let $(\alpha_{nx}, \beta_{nx})_{n=1}^\infty$ be a sequence such that $\alpha_{nx},\beta_{nx} \geq 0$ with $\lim_{n \rightarrow \infty}\alpha_{nx}+\beta_{nx}=\infty$ and $\lim_{n \rightarrow \infty}\frac{\alpha_{nx}}{\alpha_{nx}+\beta_{nx}}=\mu_x$. Then, the value function has the limit,
$$\lim_{n \to \infty}  V_x(\alpha_{nx},\beta_{nx})=\frac 1{1-\gamma_x} \max\{0,\mu_x-c\}.$$
\end{proposition}

We now introduce a preliminary structural result for an stationary optimal policy, which shows that it is optimal to forward items until some stopping time, after which we discard all items. This stopping time can be infinity, in which case we forward all items. The remark is a special case of the stopping rule shown in Theorem 5.2.2 from \cite{BeFr85}. 

\begin{remark} 
\label{re:stoppingRule}
Under any stationary optimal policy, if $U_{nx}=0$, then $U_{\ell x}=0$ for all $\ell>n$.  Thus, this optimal policy forwards all items until the stopping time $\inf \{n: U_{nx}=0\}$, and discards all subsequent items.
\end{remark}

The following lemmas state that $V_x(\alpha, \beta)$ is non-decreasing and convex as a function of $\mu(\alpha, \beta)$, holding $\alpha+\beta$ fixed. We use them below, in the proof of Theorem~\ref{th:optimalPolicy}, to show that the optimal policy for the single-category sub-problem is a threshold policy, and in the proof of Theorem~\ref{th:muStar_structure}, to show this threshold is non-decreasing.  

\begin{lemma}
\label{l:V_nonDec} 
$ \ell \mapsto V_x(\alpha+\ell, \beta-\ell)$
is non-decreasing for $-\alpha \leq \ell \leq \beta$ given any $\alpha, \beta > 0$.
\end{lemma}

\begin{lemma}
\label{l:V_convex} 
$ \ell \mapsto V_x(\alpha+\ell, \beta-\ell)$ is convex for $-\alpha \leq \ell \leq \beta$ given any $\alpha, \beta > 0$. 
\end{lemma}

Here it is the intuition. The proofs of these results use induction arguments to show that the value functions for a finite-horizon truncated version of the single-category problem is non-decreasing and convex. We then take the limit as the truncation point goes to infinity. 

Before stating the main theorem, we define a function $\mu^*(m)$ for $m>0$. Let $\mu^*(m)$ be the infimum of $\mu(\alpha, \beta)$ over all states $(\alpha, \beta)$, $\alpha, \beta >0$, with $m = \alpha+\beta$ such that it is still optimal to forward at the state. That is,
\begin{equation}
\label{eq:muStar}
\mu^*(m)=\inf \Big \{\mu(\alpha, \beta): \alpha>0, \beta>0, m=\alpha+\beta \text{ and } Q_x(\alpha, \beta, 1) > 0  \Big \}.
\end{equation}

We can think of $m$ as the effective number of observations of paper feedback, (indeed, after observing feedback on $n$ items, the corresponding value of $m$ is $m=\alpha_{nx}+\beta_{nx}=\alpha_{0x}+\beta_{0x}+n$). 
We can think of $\mu^*(m)$ as the smallest posterior mean such that we would be willing to forward. We see in the following theorem that it is optimal to forward when the posterior mean is above the threshold and discard when it is below. 

\begin{theorem}
\label{th:optimalPolicy}
  Let $(\alpha_{nx}, \beta_{nx}) = (\alpha, \beta)$ be the state of category $x$ at some step $n$ with $m=\alpha+\beta$. Then it is optimal to forward the item if $\mu(\alpha_{nx},\beta_{nx}) \ge \mu^*(\alpha_{nx}+\beta_{nx})$ and discard otherwise. In other words, the following policy is optimal:
\begin{align*}
U_{nx}=
\begin{cases}
1 & \text{if } \mu(\alpha_{nx},\beta_{nx}) \ge \mu^*(\alpha_{nx}+\beta_{nx}), \\
0 & \text{otherwise}.
\end{cases}
\end{align*}
\end{theorem}

This theorem shows that the optimal policy is a threshold policy. This provides a computational benefit, because we only need to store $\mu^*(m)$ for a one-dimensional array of possible values for $m$, rather than storing $V_x(\alpha, \beta)$ for a much larger two-dimensional array of possible values for $(\alpha,\beta)$.
The proof of the theorem is based on the monotonicity of $V_x(\alpha, \beta)$ in $\mu_x(\alpha, \beta)$, shown in Lemma \ref{l:V_nonDec}. Following the lemma, we can see that if it is optimal to forward the item at time $n$ for $\mu^*(m)$, then it is also optimal to forward for any state $(\alpha, \beta)$ with $\mu (\alpha, \beta)\geq \mu^*(m)$, since the corresponding value function is non-decreasing in $\mu(\alpha, \beta)$. 

We can compare Theorem~\ref{th:optimalPolicy} to the Gittins index policy (\cite{GiJo74}). This policy would compute the Gittins index $\nu(\alpha, \beta)$ for the (gross) value of forwarding, $\nu(\alpha, \beta)=\sup_{\tau} \frac{E[\sum_{n=1}^\tau \gamma^{n-1}Y_n | \alpha, \beta]}{E[\sum_{n=1}^\tau \gamma^{n-1}|\alpha, \beta]}$, and only forward when $\nu(\alpha, \beta) \ge c$. 
Both policies are optimal, and $\{(\alpha, \beta): \nu(\alpha, \beta) \ge c\}=\{(\alpha, \beta): \mu(\alpha, \beta) \ge \mu^*(\alpha, \beta)\}$, but $\mu^*(\cdot)$ is a function only of the effective number of samples while $\nu(\cdot, \cdot)$ depends on the full state. This offers a storage benefit, as described above. 
Also, $\mu^*(\cdot)$ can be computed by solving a single dynamic program, while computing $\nu(\cdot, \cdot)$ requires solving many dynamic programs, or doing other additional computation \citep{Gittins79, KaVe87, GiGlWe11}.

If the prior distribution on the bandit's mean reward were from a location-scale family conjugate to the bandit reward distribution, as occurs for example when the prior distribution is normal and the bandit reward is normal with known variance, then one can use linearity of the Gittins index in the posterior's location parameter (\cite{GiGlWe11} Section~7 page 192, and discussed in more detail below) to show a results similar to Theorem~\ref{th:optimalPolicy}: that the optimal policy is a threshold policy, and it is optimal to forward if and only if the posterior's location parameter is above the threshold. However, the family of Beta distributions is not a location-scale family.


Lastly, we provide structural properties of $\mu^*(m)$ in Theorem~\ref{th:muStar_structure}. 
\begin{theorem}
\label{th:muStar_structure}
$\mu^*(m)$ has the following three properties:
\begin{enumerate}
\item \indent $\mu^*(m)\le c$ for any $m>0$; 
 \item $\mu^*(m) \le \mu^*(m+1)$ for any $m>0$;
\item $\lim_{m \rightarrow \infty} \mu^*(m)=c.$
\end{enumerate}
\end{theorem}

To support the intuition behind these structural results, it is useful to decompose the net value of forwarding, $Q(\alpha,\beta,1)-Q(\alpha,\beta,0)$, into two terms:
an immediate expected reward $\mu-c$, 
and the value of information (VOI) for the observed feedback,
$\text{VOI}= \gamma_x \big[\mu V_x(\alpha+1, \beta)+(1-\mu)V_x(\alpha, \beta+1) -V_x(\alpha, \beta) \big]$, where we write $\mu=\mu(\alpha,\beta)$.
We obtain these expressions through direct examination of \eqref{eq:Q0} and \eqref{eq:Q1}.
The optimal policy forwards whenever the sum of the immediate expected reward and the value of information is non-negative.

We contrast the optimal policy with the {\it pure exploitation} policy, which ignores the value of information, and considers only the immediate expected reward.
The pure exploitation policy forwards if the immediate expected reward is non-negative, i.e., if $\mu(\alpha_{nx}, \beta_{nx}) \geq c$, and discards otherwise. This is also a threshold policy, but with a threshold of $c$.

The first part of Theorem~\ref{th:muStar_structure}, that $\mu^*(m)\le c$, shows that the optimal policy's threshold for forwarding $\mu^*(m)$ is at or below the pure exploitation policy's threshold $c$.
Thus, whenever pure exploitation forwards, the optimal policy forwards as well.  Moreover, when $\mu(\alpha,\beta)$ is in the range $[\mu^*(m),c)$, the optimal policy forwards even though the pure exploitation policy does not, exploring because the value of the feedback that will result overcomes a negative immediate expected reward.

The second part of Theorem~\ref{th:muStar_structure}, that $\mu^*(m)$ is non-decreasing, shows that this interval $[\mu^*(m),c)$ is widest, and thus also the optimal policy's willingness to forward is largest, when $m=\alpha+\beta$ is small and we have substantial uncertainty about the user's preference. As $m=\alpha+\beta$ increases, we have more feedback and less uncertainty, and the optimal policy is less willing to explore. 
A similar structural property for two-arm bandit indices was shown in Theorem 5.3.5 and Theorem 5.3.6 from \cite{BeFr85}.  These results examine how Gittins index changes as we increase the number of successes, or the number of failures, but not both.
Our result implicitly also examines how Gittins index changes as we change the effective number of measurements, but in contrast holds $\mu(\alpha, \beta)$ fixed, which requires changing both the number of successes and the number of failures.

The third part of Theorem~\ref{th:muStar_structure}, that $\lim_{m \mapsto \infty}\mu^*(m)=c$, shows that as we collect more and more feedback, and learn $\theta_x$ with greater and greater accuracy, it becomes optimal to behave like the pure exploitation policy.   This is because the pure exploitation policy is optimal when $\theta_x$ is known.

\begin{figure}[tb]
\FIGURE
  {\includegraphics*[width=3.3in, height=2.5in]{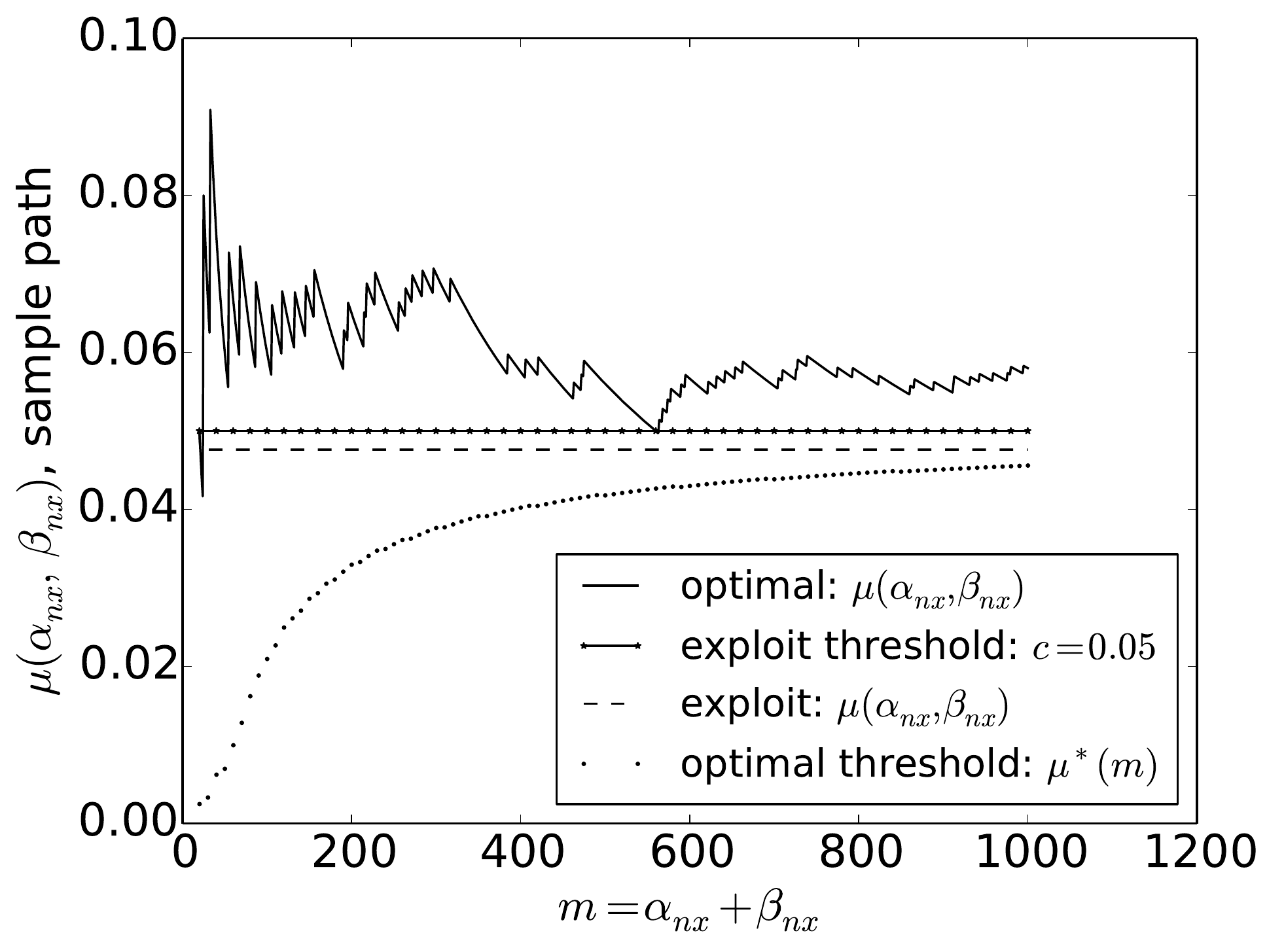}}
  {Illustration of the optimal policy's threshold, $\mu^*(m)$ (dotted line), and the pure exploitation policy's threshold, $c$ (solid line with '$*$'), with user sample paths,  $\mu(\alpha_{nx}, \beta_{nx})$, under the optimal policy (solid line) and the pure exploitation policy (dashed line), with $\alpha_{0x}=1$, $\beta_{0x}=19$, $c=0.05$ and $\gamma_x=0.999$. \label{fig:muStar_intuition}}
{}
\end{figure}

Figure~\ref{fig:muStar_intuition} shows an approximation to the optimal threshold $\mu^*(m)$, and the fixed threshold, $c$, of the pure exploitation policy in a setting with $\alpha_{0x}=1$, $\beta_{0x}=19$, $c=0.05$ and $\gamma_x=0.999$. Although $\mu^*(m)$ is non-decreasing, this approximation to $\mu^*(m)$ oscillates. This is because it is computed by taking the infimum \eqref{eq:muStar} only over values of $\alpha_{nx},\beta_{nx}$ reachable from a single $\alpha_{0x},\beta_{0x}$, rather than over all values in $(0,\infty)^2$.  
We also plot a user sample path, $\mu(\alpha_{nx}, \beta_{nx})$, under the optimal policy (solid line) and the pure exploitation policy (dashed line). 
The pure exploitation policy discards all items, because the initial value of $\mu(\alpha_{0x},\beta_{0x})$ is below $c$, but the optimal policy forwards items initially, discovers that the user's interest $\theta_x$ is larger than originally anticipated, and continues forwarding, thus earning a larger reward.

Conclusions similar to Theorem~\ref{th:muStar_structure} may be derived using standard results from the literature when the prior and posterior reside in a location-scale family, as is the case when rewards are normal with known variance, and our prior and posterior on the rewards' unknown mean is also normal.  In this case, letting $\mu$ and $\sigma^2$ be the mean and variance of the posterior, it is known (\cite{GiGlWe11} section 7 on page 192) that the Gittins index $\nu(\mu, \sigma^2)$ satisfies $\nu(\mu, \sigma^2)=\mu+\sigma \nu(0, 1)$ with $\nu(0, 1) \in (0, \infty)$. As noted above, this special structure implies a threshold policy is optimal, with threshold $\mu^*(\sigma^2)=\inf\{\mu: \nu(\mu, \sigma^2) \ge c\} = c- \sigma \nu(0, 1)$. We then immediately have $\mu^*(\sigma^2) \le c$, $\mu^*(\sigma^2)$ decreasing in $\sigma^2$, and $\lim_{\sigma^2 \rightarrow 0} \mu^*(\sigma^2)=c$. 
Thus, Theorem~\ref{th:muStar_structure} shows that properties already known for the Normal-Normal setting, and for other location-scale families, also hold in the Beta-Bernoulli setting, which is more natural for the forwarding problem.

\section{Simulation Results}
\label{sec:simulation}
In this section we present idealized simulation results for single- and multi-category problems that compare the optimal policy with the pure exploitation policy and two competitive exploration policies: upper confidence bound (UCB) and Thompson sampling.  We also evaluate the realism of these idealized simulations by comparing with more realistic trace-driven results.  Idealized simulations generate user actions and item properties from the model in Section~\ref{sec:model} using parameters estimated from historical arXiv.org data, while the trace-driven simulation uses real user histories from arXiv.org.

\subsection{Idealized Simulation} \label{sec:idealizedSim}
In this section we present Monte Carlo simulation results under five policies---optimal, pure exploitation, tuned UCB, untuned UCB, and Thompson sampling---in single-category and multi-category problems. In this setting, user histories and item properties are generated according to the assumed model from Section~\ref{sec:model} using parameters estimated from historical data from the arXiv. Because we simulate from the model, rather than using real user behavior on real items as we do in Section~\ref{sec:traceDrivenSimulation}, we call this an ``idealized simulation''.  Our numerical results demonstrate that exploration adds value in a wide variety of settings, and identify problem characteristics that determine how much value exploration adds.

In the idealized setting,
we first show the expected total reward for single-category subproblems $E^{\pi (x)} \left [\sum_{n=1}^{N_x} U_{nx}(Y_{nx}-c) \right]$. Each single-category expected total reward is estimated using $500,000$ independent simulated users, each of which is simulated for one sample path according to the model in Section~\ref{sec:model}. This simulation requires specifying four parameters: hyper-parameters $\alpha_{0x}$, $\beta_{0x}$ for the prior beta distribution, a unit forwarding cost $c$, and a discount factor $\gamma_x$. We set $\alpha_{0x}=1$ and $\beta_{0x}=19$, which are values typical among those obtained in Section~\ref{sec:traceDrivenSimulation} when fitting to historical data from arXiv.org, and which corresponds to an average user finding one out of $20$ items to be relevant. We consider a range of values for $\gamma_x$ consistent with those estimated from data in Section~\ref{sec:traceDrivenSimulation} for different arXiv categories.  For $c$, we choose a range of values near $c=0.05$. This value of $c$ is consistent with being indifferent about viewing a stream with 1 relevant item out of every $\frac 1{0.05}=20$ shown.

\begin{figure}
\FIGURE
  {\includegraphics*[scale=0.3]{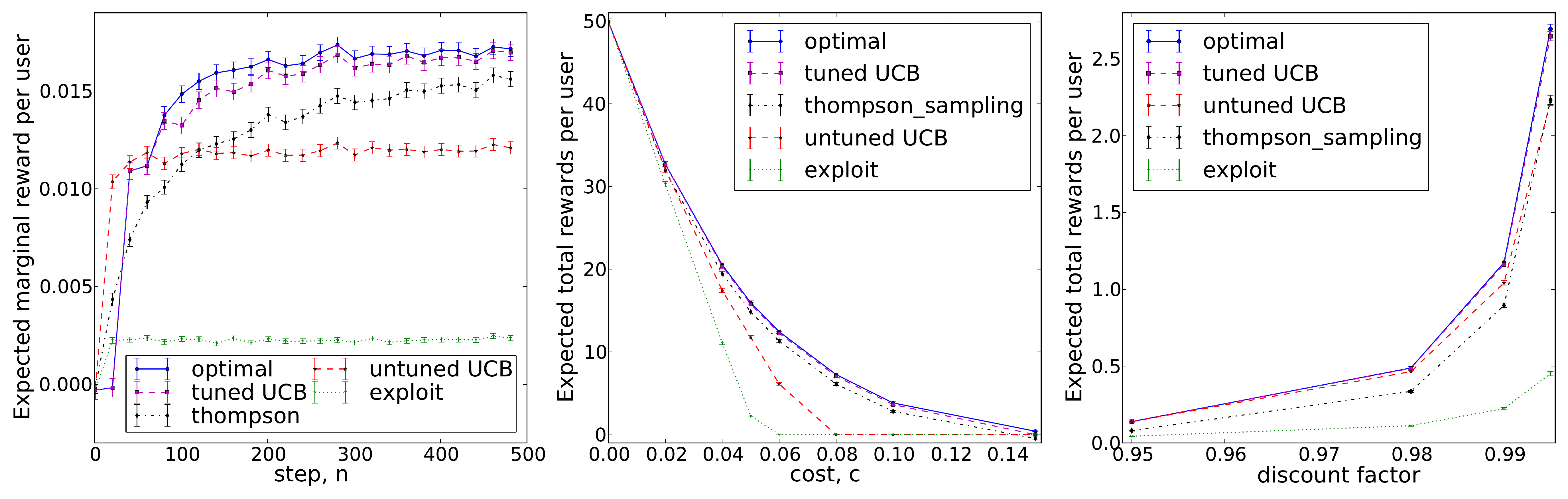}}
  {Idealized simulation results for the single-category sub-problem with $\alpha_{0x}=1$ and $\beta_{0x}=19$. The simulation compares the performance of five policies $\pi^{(x)}$: the optimal policy (denoted ``optimal''), tuned UCB, untuned UCB at $\rho=0.75$, Thompson sampling, and pure exploitation (denoted ``exploit''). Tuned UCB runs the simulation for various $\rho$ in the range $\{0.65, 0.7, 0.75, 0.8,0.85,0.9,0.95,0.99\}$ with the best results reported. Each error bar is a 95$\%$ confidence interval.  {\it (a)} The left plot shows expected marginal reward, $E^{\pi^{(x)}} \left[U_{nx}(Y_{nx}-c) \right]$, under each policy $\pi^{(x)}$ at each step $n$ with a unit forwarding cost of $c=0.05$ and a discount factor of $\gamma_x=0.999$.  {\it (b)} The middle plot shows expected total reward, $E^{\pi^{(x)}} \left[ \sum_{n=1}^{N_x} U_{nx}(Y_{nx}-c) \right]$, versus unit forward cost $c$ ranging from $0$ to $0.15$ with a discount factor of $\gamma_x=0.999$.  {\it (c)} The right plot shows expected total reward versus discount factor, $\gamma_x$, ranging from $0.95$ to $0.995$ with a unit forwarding cost of $c=0.05$. \label{fig:idealizedSimGraphs} }
{}
\end{figure}

The results of this simulation for the single-category problems under the five policies, optimal, pure exploitation, tuned UCB, untuned UCB, and Thompson sampling, are plotted in Figure \ref{fig:idealizedSimGraphs}, with 95$\%$ confidence intervals shown as error bars. UCB and Thompson sampling are two heuristic exploration policies based on common approaches to exploration vs. exploitation in the broader literature \citep{LaRo85, AuCeFi02, Kaufmann2012, Thompson1933, AgGo2011, ChapelleLi2011, RussoVanRoy2014}. At each step the UCB policy computes the $\rho$-quantile, $Q(\rho, \theta_x)$, associated with the posterior distribution of $\theta_x$ and forwards the item if $Q(\rho, \theta_x) \ge c$, while the Thompson sampling policy draws a sample, $\hat{\theta}_x$, from the posterior and forwards if $\hat{\theta}_x \ge c$.   The ``tuned UCB'' policy refers to the UCB policy where $\rho$ is tuned for each given set of problem parameters ($\gamma_x$, $\alpha_{0x}$, $\beta_{0x}$, $c$: $x=1,...,k$) by using simulation to try several $\rho$ in the range $\{0.65, 0.7, 0.75, 0.8, 0.85, 0.9, 0.95, 0.99\}$ and using the value that performs best while in the ``untuned UCB'' policy $\rho=0.75$ is chosen  and fixed for all $c$ and $\gamma_x$ in Figure~\ref{fig:idealizedSimGraphs} since $\rho=0.75$ is best for smaller unit costs and discount factors. 

Figure \ref{fig:idealizedSimGraphs}(a) shows the expected marginal reward $E^{\pi(x)} \left[ U_{nx}(Y_{nx}-c)\right]$ at each time step $n \in \{1, .., 500\}$ with a cost of $c=0.05$ and a discount factor of $\gamma_x=0.999$. In this setting, tuned UCB chooses $\rho=0.95$ while untuned UCB sets $\rho=0.75$. In the earlier steps, the expected marginal reward of the optimal policy and tuned UCB suffers from extensive exploration compared to pure exploitation, Thompson sampling, and untuned UCB. But as we start collecting information about users over the next 50 steps, the expected marginal reward of the optimal policy and tuned UCB quickly recovers from the previous loss and rapidly surpasses that of the other policies, eventually stabilizing as feedback becomes abundant, and the optimal policy stops exploring and exploits the information it has gained.

\begin{figure}
\FIGURE
  {\includegraphics*[scale=0.4]{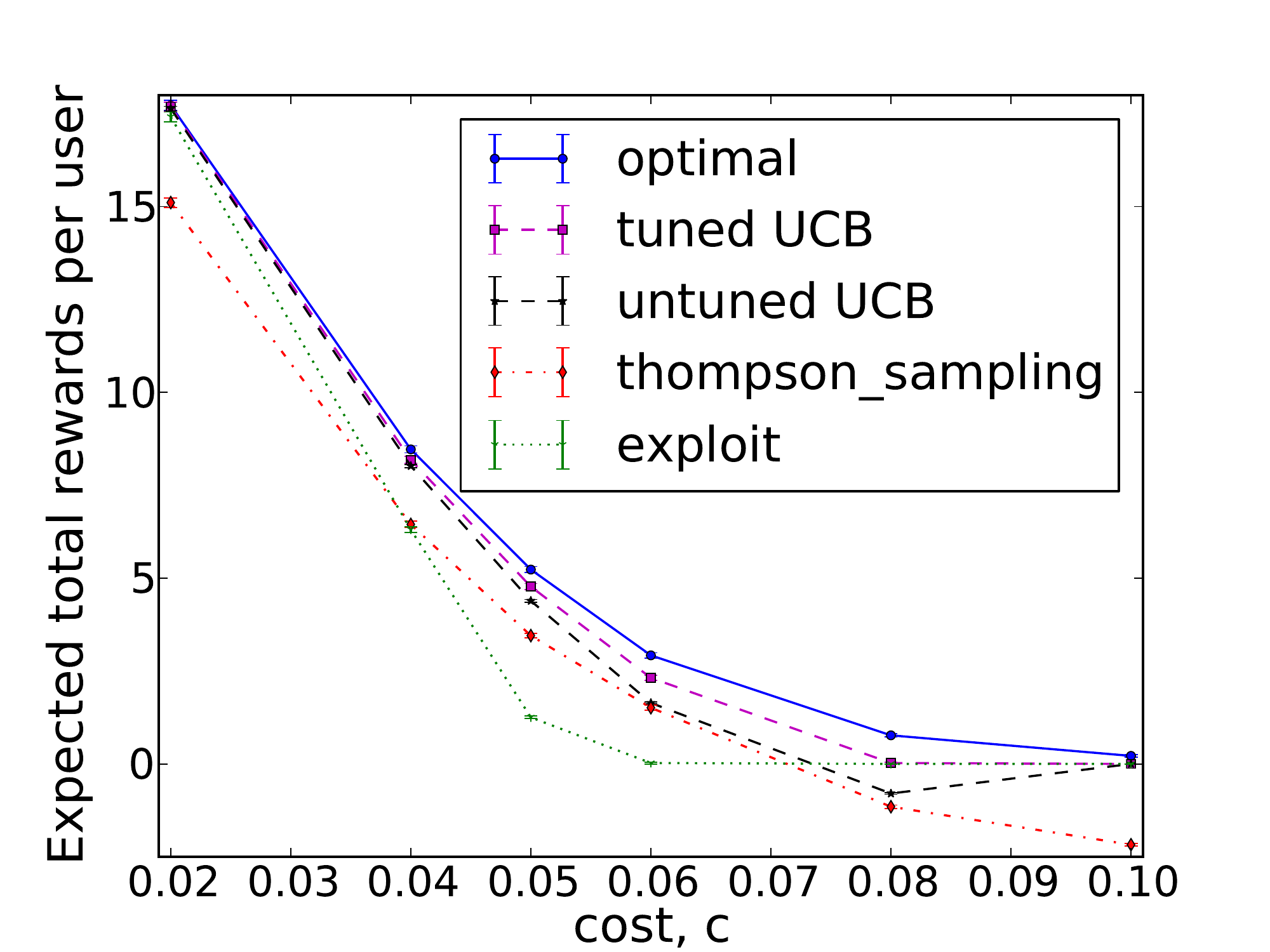}}
  {Idealized simulation results in a multi-category problem, with a mixture of 20 categories at $(\alpha_{0x}, \beta_{0x}, \gamma_{0x})=(1,19,0.95)$ and one category at $(\alpha_{0x}, \beta_{0x}, \gamma_{0x})=(1,19,0.995)$. The simulation compares the performance of five policies $\pi$: the optimal policy, tuned UCB, untuned UCB with $\rho=0.85$, pure exploitation, and Thompson sampling. Each error bar is a 95$\%$ confidence interval.  The plot shows expected total reward, $E^{\pi} \left[ \sum_{n=1}^{N} U_{n}(Y_{n}-c) \right]$, versus unit forwarding cost $c$ ranging from $0.02$ to $0.1$. In tuned UCB, simulations are run for a range of $\rho$-quantiles, $\{0.65, 0.7, 0.75,0.8,0.85,0.9,0.95,0.99\}$, with the best expected total rewards reported in the figure for each cost, $c$. Untuned UCB sets $\rho=0.85$ since it performs the best in the category with $\gamma_{0x}=0.995$. \label{fig:idealizedSim_multiCats} }
{}
\end{figure}

In addition to the expected marginal reward, Figure \ref{fig:idealizedSimGraphs}(b) and Figure \ref{fig:idealizedSimGraphs}(c) show sensitivity plots of the expected total reward. Figure \ref{fig:idealizedSimGraphs}(b) shows the expected total reward against different unit costs $c$ for a fixed discount factor $\gamma_x=0.999$, while Figure \ref{fig:idealizedSimGraphs}(c) shows the expected total reward per user against different discount factors $\gamma_x$ with a fixed unit cost $c=0.05$. In the single-category subproblems, we observe that the optimal policy is almost identical to tuned UCB, and performs statistically better than Thompson sampling, untuned UCB, and pure exploitation. The difference between the optimal policy and tuned UCB becomes statistically significant in multi-category problems, as illustrated in Figure~\ref{fig:idealizedSim_multiCats}, because different optimal tunings of $\rho$ are required for categories with different discounts. 
In Figure~\ref{fig:idealizedSimGraphs}(b), the expected total reward decreases as cost increases for all five policies, which is intuitive because each forwarded item provides a smaller (net) reward when costs are higher. The optimal policy and tuned UCB performs at least as well as the other three policies for all cost values. At the two extremes when the cost is close to either 0 or 1, all five policies make the same forwarding decision, i.e. to forward all available items when the cost is near $0$, or forward none when the cost is near 1. This explains why the five policies have the same expected total reward near the two extremes. As we move away from the extremes, the optimal policy provides a greater competitive advantage, by optimally balancing exploration and exploitation. It outperforms pure exploitation the most when the unit cost is near $c=0.05$, because exploration provides the most benefit in the most ambiguous case when $c$ is close to the expected value of $\theta_x$, and our prior has mean $\mu(\alpha_{0x}, \beta_{0x})=\frac{1}{1+19}=0.05$. Compared to tuned UCB, fixing $\rho=0.75$ for all $c$ in untuned UCB deteriorates the performance.

In Figure \ref{fig:idealizedSimGraphs}(c), we observe that the expected total reward increases with the discount factor, $\gamma_x$, for all policies, but the optimal policy and tuned UCB have a much steeper positive slope as the discount factor approaches $1$. This suggests that the optimal policy benefits the most from exploration when the discount factor is large, since the system has a longer time horizon to learn users' preferences and recover any losses suffered at the beginning when learning. 

In Figure~\ref{fig:idealizedSim_multiCats} we compare the performance of the optimal policy with the other four policies for the idealized simulation in a multi-category problem. The constructed multi-category problem consists of 20 categories at $(\alpha_{0x}, \beta_{0x}, \gamma_x)=(1,19,0.95)$ and one category at $(\alpha_{0x}, \beta_{0x}, \gamma_x)=(1,19,0.995)$. Untuned UCB sets $\rho=0.85$ since it is best when $\gamma_{0x}=0.995$. Figure~\ref{fig:idealizedSim_multiCats} shows that tuned UCB, with the flexibility of tuning $\rho$, is comparable to the optimal policy when the unit forwarding cost is near 0, but its difference from the optimal policy enlarges and becomes statistically significant as the unit forward cost increases. Tuning UCB improves its performance over untuned UCB, but does not bring it as close to optimal as it did in the single-category case. The reason is that the right balance between exploration vs. exploitation differs across item category, preventing UCB from achieving this balance with a single value of $\rho$. 

One could imagine UCB with a more elaborate tuning method, in which $\rho$ is allowed to differ across item category, and is tuned separately using simulation optimization for each user and category based on $\alpha_{0x}$, $\beta_{0x}$, $\gamma_x$, $c$, but using this heuristic approach would not offer any significant computational advantage over simply using the optimal policy.

\subsection{Trace-driven Simulation}\label{sec:traceDrivenSimulation}
In this section, we present trace-driven simulation results using the web server logfile from arXiv.org for both the optimal policy and the pure exploitation policy. 
Instead of using simulated users and items as in Section \ref{sec:idealizedSim}, we use real historical user interactions and items extracted from this logfile.  All identifiable information from users is hashed to protect user privacy. In contrast with the focus on various policies in Section \ref{sec:idealizedSim}, we choose only the optimal policy and the pure exploitation policy in this section to test strategies at the two ends of the spectrum of performance, to understand better how exploration adds value in a more realistic setting and how predictions from our model match reality. 

\begin{figure}
\FIGURE
{\includegraphics*[scale=0.5]{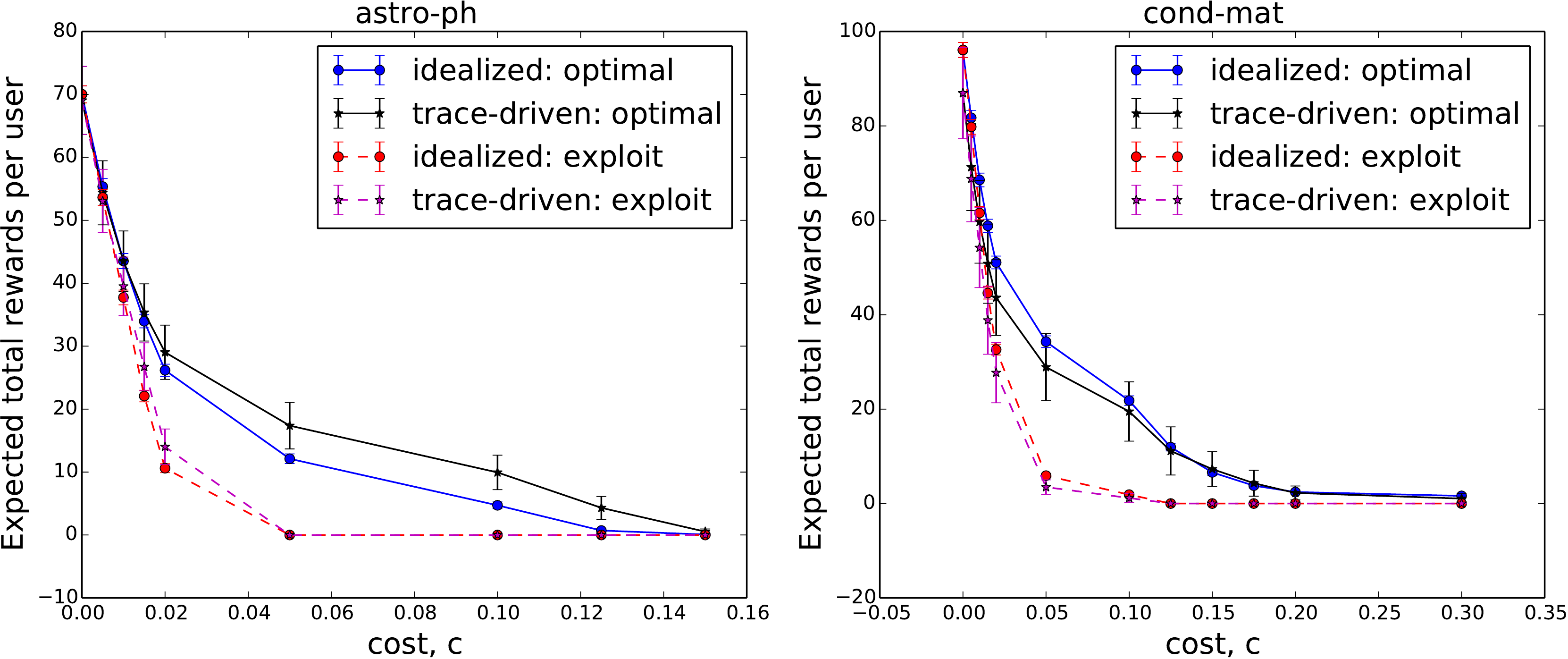}}
{This figure plots expected cumulative reward against the unit cost of forwarding $c$ under the trace-driven simulation and the idealized Monte-Carlo simulation for (a) astrophysics (astro-ph, left), and (b) condensed matter (cond-mat, right).  Results are presented for both the optimal policy and the pure exploitation policy. 95$\%$ confidence intervals are shown as error bars.\label{fig:traceDrivenSimGraph}}
{}
\end{figure}

ArXiv items are categorized by their authors into one of 18 subjects (astrophysics, condensed matter, computer science, ...) and then again into one of several categories within the subject (for example, the condensed subject has nine categories: Disordered Systems and Neural Networks, Materials Science, Mesoscale and Nanoscale Physics, ...). In addition to this ``primary'' subject/category, authors may optionally provide one or more secondary subject/category labels. 

Each subject has associated ``new'' and ``recent'' webpages, which show items submitted during the previous day and week, respectively. Some large subjects, e.g., astrophysics and condensed matter, have over 100 items submitted per day. In our trace-driven simulation, we take items submitted to a single subject as our item stream, and use the author-provided primary category within this subject as our item category. We use this simple pre-existing categorization method to focus attention on the exploration vs. exploitation tradeoff rather than the categorization scheme, although one could easily use categories learned automatically from item content and/or historical co-access data \citep{IR_manning2008}.

In our experiments, we consider two separate item streams, astrophysics and condensed matter, and look at items submitted in 2009 and 2010. In each of these subjects, we consider an item to be presented to the user if he/she visited the subject's new or recent page during the period it was posted there. If the user clicked on the link to the full text or abstract from the new or recent page during this period, we consider the item to be relevant. If he/she did not click on this link, we consider the item to be irrelevant. If the user did not visit the new or recent page during the period the item was posted there, we consider the item's relevance to be unobserved. We then identify those users who visited the subject's new and recent pages a moderate number of items over the time period of interest (2009-2010), removing those who visited too infrequently (less than 30 visits) as not providing useful data, and those who visited too often (more than 510 visits) as likely robots. 

In the next step, we randomly assign each extracted user into one of two groups: training and testing. The training users are used to estimate hyper--parameters $\alpha_{0x}$, $\beta_{0x}$, $\gamma_x$ for each category $x$, and then the simulation is executed among the testing users. For each training user $u$, we first estimate his/her $\theta_{u, x}$ associated with each category $x$ (we have added $u$ to the subscript for $\theta_x$, and below for $N_x$, to emphasize the dependence on the user) by calculating his/her click-through-rate on the items presented from that category. We then estimate $\alpha_{0, x}$ and $\beta_{0, x}$ for each category $x$ by fitting a beta distribution to the histogram of estimated $\theta_{u, x}$ over all training users $u$. Similarly, to estimate $\gamma_x$ for each category $x$, we first count the total number $N_{u, x}$ of items submitted to the category over each user $u$'s lifetime in the system, and fit a geometric distribution with parameter $1-\gamma_x$ to the histogram of $N_{u,x}$ among the training users. 

After the parameter estimation, we then perform the trace-driven simulation for each testing user as follows. Iterating through those dates on which the user visited the category's new or recent page, we take all the items that were shown on those pages, and decide sequentially whether to forward that the item to the user according to one of the two policies. If the user clicked, according to the historical data, on the link to the full text or the abstract of the forwarded item over the period when the item was posted on the new or recent page, we conclude that the user found the item to be relevant; otherwise, if the user did not click, we conclude that the item was irrelevant to the user. Feedback is collected immediately and given to the policy, and the process is repeated for the next item, until the user leaves the system.

For a side-by-side comparison, we also perform an idealized simulation for the multi-category problem to check how trace-driven results differ from the idealized results. This idealized simulation follows the same framework as in Section \ref{sec:idealizedSim} but uses multiple categories, and uses the same values for $\alpha_{0x}$, $\beta_{0x}$, and $\gamma_x$ as the trace-driven simulation.

Figure \ref{fig:traceDrivenSimGraph} shows the expected total reward per user versus the unit forwarding cost $c$ under the two simulations, with error bars representing 95$\%$ confidence intervals. 
The patterns observed are similar to those in Figure \ref{fig:idealizedSimGraphs}(b): the two policies perform similarly when the unit cost is near $0$ or 1; and the optimal policy substantially outperforms pure exploitation when $c$ is away from the two extremes. 

For each policy, the trace-driven results are close to the idealized simulation results, with much of the discrepancy explained by sampling error.
However, some of the discrepancy is likely due to violations of our modeling assumptions by the historical data.
In our model, we make the four main modeling assumptions: (1) arrivals of items follow a Poisson process and the lifetime of a user in the system is exponentially distributed, so that the total number of items in a category viewed by each user, $N_x$, follows a geometric distribution; (2) the prior distribution on $\theta_x$ in each category follows a Beta distribution;  (3) the $\theta_x$ are independent across categories; (4) the number of items in the category $x$ viewed by the user, $N_x$, and probability of relevance of an item from category $x$ to the user, $\theta_x$,  are independent. We performed some empirical checks on the historical data to validate these assumptions. Assumptions (1) and (2) seem to be met reasonably well, but we saw some violations of (3) and (4). We leave extensions of our model to incorporate more general assumptions to future work.

Despite these violations of our modeling assumptions, we see that simulation model matches well with the behavior of the more realistic trace-driven simulation, and that the policy calculated to be optimal in our model provides similar improvements over pure exploitation in both trace-driven and idealized simulations. 

\section{Conclusion}
We formulated a personalized information filtering problem as a Bayesian sequential decision-making problem and provided the optimal policy for adaptively forwarding items from multiple categories to maximize the expected total reward via balancing exploration against exploitation. With an independence property, the original $k$-category problem can be decomposed to an aggregation of $k$ single-category subproblems, avoiding the so-called ``curse of dimensionality". Moreover, we show that the optimal policy for each sub-problem is a threshold policy, and that this threshold has intuitive structural properties. To compute this threshold, we provide an approximation to the optimal policy with rigorous error bounds, whose error converges geometrically to  0 as truncation increases. Lastly, results from both idealized Monte Carlo and trace-driven simulations show that this optimal policy provides value in practice. 

\section*{Appendix A. Computation of the single-category value function}
We cannot compute the single-category value function $V_x(\alpha, \beta)$ exactly through Bellman's recursion because storing $V_x(\alpha, \beta)$ for all possible values of $\alpha$ and $\beta$ would require infinite storage.  In this section we describe a method for computing an approximation with rigorous error bounds for $V_x(\alpha,\beta)$, from which an approximation to the optimal single-category policy can be computed.

The following lemma follows directly from Proposition~\ref{prop:lowerBound} and~\ref{prop:upperBound}, where we have loosened the upper bound for easier computation using the inequality $E[ \max\{\theta_x-c, 0\}]\leq 1$.
\begin{lemma} \label{lemma:bounds}
  For all $\alpha,\beta>0$, we have 
  $\frac1 {1-\gamma_x} \max\{0,\mu(\alpha,\beta)-c\}
  \le V_x(\alpha,\beta) \le 
  \frac1{1-\gamma_x}$.
\end{lemma}

We now introduce an algorithm, Algorithm~\ref{alg1}, that defines and computes two quantities, $V_x^U(\alpha,\beta; M)$ and $V_x^L(\alpha,\beta;M)$, which bound the value function $V_x(\alpha,\beta)$ above and below, as stated below in Lemma~\ref{l:approximationBounds}. $V^U_x(\alpha, \beta;M)$ and $V_x^L(\alpha, \beta;M)$ are value functions for a finite horizon truncated version of the problem, where after being presented with M items we are given a terminal reward, which is the lower bound $\frac1{1-\gamma_x}\max\{0, \mu(\alpha_{Mx}, \beta_{Mx})\}$ from Lemma~\ref{lemma:bounds} when calculating $V_x^L$, and the upper bound  when calculating $V_x^U$. 

\begin{algorithm}[!h]
\caption{Computation of $V_x^L(\alpha, \beta; M)$ and $V_x^U(\alpha, \beta; M)$}
\label{alg1}
\begin{algorithmic}
\REQUIRE $\alpha_{0x}$, $\beta_{0x}$, $\gamma_x$, $c$, and $M$
\FOR{$i=0,...,M$}
	\STATE{Let $\alpha=\alpha_{0x}+i$, $\beta=\beta_{0x}+M-i$.}
	\STATE{Let $V^L_x(\alpha, \beta; M)=\frac 1{1-\gamma_x} \max \left\{0, \frac{\alpha}{\alpha+\beta}-c\right\}$ and 
	$V^U_x(\alpha, \beta; M)=\frac 1{1-\gamma_x}$.}
\ENDFOR
 \FOR{$\ell=M-1, M-2, ..., 0$}
	\FOR{$i=0, ..., \ell$}
	\STATE{Let $\alpha=\alpha_{0x}+i$, $\beta=\beta_{0x}+\ell-i$.}
	\STATE{Let $V^L_x(\alpha, \beta; M)=\max \left\{0, \frac{\alpha}{\alpha+\beta}-c+\gamma_x \left[\frac{\alpha}{\alpha+\beta}V^L_x(\alpha+1, \beta; M)+\frac{\beta}{\alpha+\beta}V^L_x(\alpha, \beta+1; M)\right] \right\}$}.
	\STATE{Let $V^U_x(\alpha, \beta; M)=\max \left\{0, \frac{\alpha}{\alpha+\beta}-c+\gamma_x \left[\frac{\alpha}{\alpha+\beta}V^U_x(\alpha+1, \beta; M)+\frac{\beta}{\alpha+\beta}V^U_x(\alpha, \beta+1; M)\right] \right\}$}.
	\ENDFOR
  \ENDFOR
\end{algorithmic}
\end{algorithm}

\begin{lemma}
\label{l:approximationBounds}
  For each $M\ge 0$, $0\le \ell \le M$, and $0 \le i \le \ell$ with $\alpha = \alpha_{0x} + i$ and $\beta = \beta_{0x} + \ell - i$,
  we have
$V^L_x(\alpha,\beta; M) \le V_x(\alpha,\beta) \le V^U_x(\alpha,\beta; M)$.
\end{lemma}

These two bounds provide a computable approximation to the value function, and thus to the optimal policy. The next lemma shows that the gap between these two bounds converges to 0 as the termination step $M$ approaches infinity. 

\begin{lemma}
\label{l:gap_convergence}
  For each $M\ge 0$, $0\le \ell \le M$, $0 \le i \le \ell$ with $\alpha = \alpha_{0x} + i$ and $\beta = \beta_{0x} + \ell - i$, we have
  $V^U_x(\alpha, \beta; M) - V^L_x(\alpha, \beta; M) \le \gamma_x^{M-\ell} / (1-\gamma_x)$. 
In particular, 
$\lim_{M\rightarrow \infty} V^U_x(\alpha, \beta; M) - V^L_x(\alpha, \beta; M) = 0$.
\end{lemma}

\section*{Appendix B. Mathematical Proofs}
\subsection*{Proof of Theorem~\ref{t:decomposition}}
 Let $V$ be equal to the value of \eqref{eq:overall-problem}, and for each $x$, let $V_x$ be equal to the value of \eqref{eq:single-category-problem}. (We have dropped the $\alpha$ and $\beta$ from the notation $V$ and $V_x$ in this proof.) By construction, on each sample path,
  \begin{equation*}
    \sum_{n=1}^N U_n (Y_n-c) = \sum_{x=1}^k \sum_{\ell=1}^{N_x} U_{\ell x}(Y_{\ell x}-c).
  \end{equation*}
  Thus, for each $\pi\in\Pi$,
  \begin{equation}
    \label{eq:proof-1a}
    E^\pi\left[ \sum_{n=1}^N U_n (Y_n-c)\right]
    = E^\pi\left[ \sum_{x=1}^k \sum_{\ell=1}^{N_x} U_{\ell x}(Y_{\ell x}-c) \right]
    = \sum_{x=1}^k E^\pi\left[ \sum_{\ell=1}^{N_x} U_{\ell x}(Y_{\ell x}-c) \right],
  \end{equation}
  and we have
  \begin{equation*}
    V = \sup_{\pi\in\Pi} E^\pi\left[ \sum_{n=1}^N U_n (Y_n-c)\right]
      = \sup_{\pi\in\Pi} \sum_{x=1}^k E^\pi\left[ \sum_{\ell=1}^{N_x} U_{\ell x}(Y_{\ell x}-c) \right]
      \le \sum_{x=1}^k \sup_{\pi\in\Pi} E^\pi\left[ \sum_{\ell=1}^{N_x} U_{\ell x}(Y_{\ell x}-c) \right],
    \end{equation*}
    where the last inequality follows from the fact that the right-hand side potentially allows a different $\pi$ to attain the supremum for each $x$.

Now fix an $x$ and consider the term
$\sup_{\pi\in\Pi} E^\pi\left[ \sum_{\ell =1}^{N_x} U_{\ell x}(Y_{\ell x}-c) \right]$.
This is equal to
$V_x=\sup_{\pi^{(x)}\in\Pi^{(x)}} E^{\pi^{(x)}}\left[ \sum_{\ell=1}^{N_x} U_{\ell x}(Y_{\ell x}-c) \right]$
because the conditional distribution of $(Y_{ix}, N_x: i\ge \ell)$
given the history available to any policy $\pi\in\Pi$ when making the decision $U_{\ell x}$
(this history is $(X_i, U_i, U_iY_i : i<n_{\ell x})$ and $X_{n_{\ell x}}=x$)
depends only upon the history available to a policy $\pi^{(x)}\in\Pi^{(x)}$, which is
$(U_{ix},U_{ix}Y_{ix}: i<\ell)$.

Thus, we have
$\sup_{\pi\in\Pi} E^\pi\left[ \sum_{\ell=1}^{N_x} U_{\ell x}(Y_{\ell x}-c) \right] =V_x$,
and
\begin{equation}
  \label{eq:proof-1b}
  V\le \sum_{x=1}^k V_x.
\end{equation}

We now show the opposite inequality, and the additional claim about constructing an optimal policy for \eqref{eq:overall-problem}.  Let $\pi^{(x),*}$ and $\pi^*$ be as described in the statement of the theorem.  Then, for each $x$,
\begin{equation*}
E^{\pi^*}\left[ \sum_{\ell=1}^{N_x} U_{\ell x}(Y_{\ell x}-c) \right] =
E^{\pi^{(x),*}}\left[ \sum_{\ell=1}^{N_x} U_{\ell x}(Y_{\ell x}-c) \right] = V_x.
\end{equation*}
Summing over $x$ and using \eqref{eq:proof-1a}, we have that the value of $\pi^*$ in the original problem is
\begin{equation}
  \label{eq:proof-1c}
  E^{\pi^*}\left[ \sum_{n=1}^N U_n(Y_n-c) \right]
  = \sum_{x=1}^k E^{\pi^*}\left[ \sum_{\ell=1}^{N_x} U_{\ell x}(Y_{\ell x}-c) \right]
  = \sum_{x=1}^k V_x.
\end{equation}

Since $\pi^*\in\Pi$, and $V$ is defined by taking the supremum over all policies $\pi$, we also have
\begin{equation}
\label{eq:proof-1d}
E^{\pi^*}\left[ \sum_{n=1}^N U_n(Y_n-c) \right] \le V.
\end{equation}
Combining
\eqref{eq:proof-1b},
\eqref{eq:proof-1c}
and \eqref{eq:proof-1d}
we have,
\begin{equation*}
\sum_{x=1}^k V_x
= E^{\pi^*}\left[ \sum_{n=1}^N U_n(Y_n-c) \right]
\le V
\le \sum_{x=1}^k V_x,
\end{equation*}
which implies both that $V=\sum_{x=1}^k V_x$, and that $\pi^*$ attains the supremum defining $V$, as claimed.

\subsection*{Proof of Remark~\ref{r:geometric}}
Fix any $n \in \{0, 1, 2, \ldots\}$, 
\begin{align*}
  P(N_x = n)
  &= \sum_{m \ge n} P(N = m) {m \choose n} p_z^n (1-p_z)^{m-n}\\
  &= \sum_{m \ge n} \gamma^m (1-\gamma){m \choose n} p_z^n (1-p_z)^{m-n}\\
  &= (1-\gamma) (p_z \gamma)^n \sum_{m \ge n} {m \choose n} [(1-p_z)\gamma]^{m-1}\\
  &=\frac{(1-\gamma) (p_z \gamma)^n}{[1-(1-p_z)\gamma]^{n+1}}
  =\gamma_x^n (1-\gamma_x).
\end{align*}
The third equality use the formula that $\frac{1}{(1-x)^s} = \sum_{k=0}^\infty {s+k-1 \choose k} x^k = \sum_{k=0}^\infty {s+k-1 \choose s-1} x^k$. With that, we show the claim. 

\subsection*{Proof of Lemma~\ref{l:discount}}
\begin{equation}
  \label{eq:proof-2a}
    E^{\pi^{(x)}}\left[\sum_{n=1}^{N_x} U_{nx}(Y_{nx}-c)\right]
  = E^{\pi^{(x)}}\left[\sum_{n=1}^\infty \indicator{n \le N_x} U_{nx}(Y_{nx}-c)\right]
  = \sum_{n=1}^\infty E^{\pi^{(x)}}\left[\indicator{n \le N_x} U_{nx}(Y_{nx}-c)\right],
\end{equation}
where the last equality is justified by Fubini's theorem and 
$E^{\pi(x)}\left[\sum^\infty_{n=1}|\indicator{n \le N_x} U_{nx}(Y_{nx}-c)| \right]<\infty$.

Considering one of these terms, for any fixed $n$, we have that 
$E^{\pi^{(x)}}\left[\indicator{n \le N_x} U_{nx}(Y_{nx}-c)\right]$
is equal to
\begin{align*}
  E^{\pi^{(x)}}\left[ E^{\pi^{(x)}}\left[\indicator{n \le N_x} U_{nx}(Y_{nx}-c) \mid U_{nx},Y_{nx} \right]\right]
  &= E^{\pi^{(x)}}\left[ P^{\pi^{(x)}}\left(n \le N_x \mid U_{nx},Y_{nx} \right) U_{nx}(Y_{nx}-c) \right]\\
  &= E^{\pi^{(x)}}\left[ \gamma_x^n U_{nx}(Y_{nx}-c) \right].
\end{align*}
Plugging this expression into \eqref{eq:proof-2a} and applying Fubini's theorem again shows
\begin{align*}
    E^{\pi^{(x)}}\left[\sum_{n=1}^{N_x} U_{nx}(Y_{nx}-c)\right]
    &= \sum_{n=1}^\infty E^{\pi^{(x)}}\left[ \gamma_x^n U_{nx}(Y_{nx}-c) \right] = \gamma_x E^{\pi^{(x)}}\left[ \sum_{n=1}^\infty \gamma_x^{n-1} U_{nx}(Y_{nx}-c) \right],
\end{align*}
which is the claimed expression.

\subsection*{Proof of Remark~\ref{re:notForward}}
When it is optimal to discard at $(\alpha, \beta)$, we have $V_x(\alpha, \beta)=Q_x(\alpha, \beta, 0)=\gamma_x \cdot V_x(\alpha, \beta)$ for some $0<\gamma_x<1$. This implies that $V_x(\alpha, \beta) = Q_x(\alpha, \beta, 0)=0$.  

We now show the other direction. If $V_x(\alpha, \beta) = 0$, then $Q_x(\alpha, \beta, 0) = \gamma_x V_x(\alpha, \beta) = 0$ and furthermore $Q_x(\alpha, \beta, 1) \leq V_x(\alpha, \beta) = Q_x(\alpha, \beta, 0) = 0$. 

\subsection*{Proof of Proposition~\ref{prop:lowerBound}}
To show the lower bound, we consider a policy that ignores learning and makes decisions based on the conditional expected reward, $ E\left[Y|\theta \sim \Bta(\alpha, \beta)\right]=\mu(\alpha,\beta)=\frac{\alpha}{\alpha+\beta}$. Let us define the policy, $\pi_L^{(x)}$, such that for all $n\geq 1$:
\begin{align*}
U_{nx}=
\begin{cases}
1 & \text{if } \mu(\alpha, \beta) \geq c, \\
0 & \text{otherwise}.
\end{cases}
\end{align*}
Then, this policy has value:
\begin{align*}
V_x^{L}(\alpha,\beta)
&=E^{\pi_L^{(x)}}\left[\sum_{n=1}^\infty \gamma_x^{n-1} U_{n x} (Y_{n x}-c) \mid \theta_x \sim \mathrm{Beta}(\alpha,\beta)\right] \\
&=\sum^\infty_{n=1}\gamma_x^{n-1} ~ \max \Big \{0, \mu(\alpha,\beta)-c \Big \} =\frac 1{1-\gamma_x}\max \big \{0, ~\mu(\alpha,\beta)-c \big \} \\
& \leq \sup_{\pi^{(x)} \in \Pi^{(x)}} E^{\pi^{(x)}}\left[\sum_{n=1}^\infty \gamma_x^{n-1} U_{n x} (Y_{n x}-c) \mid \theta_x \sim \mathrm{Beta}(\alpha,\beta)\right] =V_x(\alpha,\beta) 
\end{align*}

\subsection*{Proof of Proposition~\ref{prop:upperBound}}
We have
\begin{align*}
V_x(\alpha,\beta)
&=\sup_{\pi^{(x)} \in \Pi^{(x)}} E^{\pi^{(x)}} \left [ \sum^\infty_{n=1} \gamma_x^{n-1} (Y_{n x}-c)U_{n x} ~|~ \alpha_{0x}=\alpha, \beta_{0x}=\beta \right ] \\
&=\sup_{\pi^{(x)} \in \Pi^{(x)}} \sum^\infty_{n=1}\gamma_x^{n-1} E^{\pi^{(x)}} \Big[(Y_{n x}-c)U_{n x} ~|~ \alpha_{0x}=\alpha, \beta_{0x}=\beta \Big] \\
&=\sup_{\pi^{(x)} \in \Pi^{(x)}} \sum^\infty_{n=1}\gamma_x^{n-1} E^{\pi^{(x)}} \Big[ E^{\pi^{(x)}} \big[(Y_{n x}-c)U_{n x} ~|~ H_{n-1, x}, \theta_x, \alpha_{0x}=\alpha, \beta_{0x}=\beta \big] ~|~ \alpha_{0x}=\alpha, \beta_{0x}=\beta \Big] \\
&=\sup_{\pi \in \Pi^{(x)}} \sum^\infty_{n=1}\gamma_x^{n-1} E^{\pi^{(x)}} \Big [(\theta_x-c)U_{nx} ~|~ \alpha_{0x}=\alpha, \beta_{0x}=\beta \Big ] \\
&\leq \sum^\infty_{n=1}\gamma_x^{n-1} E \Big [\max \left \{0, \theta_x-c \right\} ~|~ \alpha_{0x}=\alpha, \beta_{0x}=\beta \Big ] \\
&=\frac 1{1-\gamma_x} E \Big [\max \left \{0, \theta_x-c \right\} ~|~ \alpha_{0x}=\alpha, \beta_{0x}=\beta \Big ]=V_x^U(\alpha, \beta).
\end{align*}

The first equality is the definition of the value function in equation~\eqref{eq:valueFunction},
while the second equality is due to Fubini's Theorem since 
$E^{\pi(x)} \left[ \sum^\infty_{n=1}|\gamma_x^{n-1}(Y_{nx}-c)U_{nx}| \big | \alpha_{0x}=\alpha \right]$.
Then applying the tower property of conditional expectation, we derive the third equation. Because the expected value of $Y_{n x}$ conditioned on $\theta_x$ is $\theta_x$ and $U_{n x}$ only depends on $H_{n-1,x}$, the inner conditional expectation reduces to $(\theta_x-c)U_{n x}$, which is smaller than $\max\{0, \theta_x-c\}$. So the inequality holds in the fifth line and we derive an upper bound for the value function. 

\subsection*{Proof of Proposition~\ref{prop:convergence}}

\noindent First, we show a lower bound on $\liminf_{n\to \infty} V_x(\alpha_{nx}, \beta_{nx})$:
\begin{align*}
\liminf_{n \to \infty} V_x(\alpha_{n x},\beta_{n x}) 
&\ge \liminf_{n \to \infty} \left[\frac 1{1-\gamma_x}\max \left\{0, \frac{\alpha_{n x}}{\alpha_{n x}+\beta_{n x}}-c \right \} \right] \\
&= \lim_{n \to \infty}  \left[\frac 1{1-\gamma_x}\max \left\{0, \frac{\alpha_{n x}}{\alpha_{n x}+\beta_{n x}}-c \right \} \right] \\
&= \frac1{1-\gamma_x} \max \left\{0,\lim_{n \rightarrow \infty} \frac {\alpha_{n x}}{\alpha_{n x}+\beta_{n x}}-c \right\} \\
&= \frac1{1-\gamma_x} \max\{0,\mu_x-c\}
\end{align*} 

The first inequality uses Proposition~\ref{prop:lowerBound}, and the rest follows because the limit of $\frac 1{1-\gamma_x}\max \left\{0, \frac{\alpha_{n x}}{\alpha_{n x}+\beta_{n x}}-c \right \}$ exists. Similarly, using Proposition~\ref{prop:upperBound}, we show an upper bound on $\limsup_{n \to \infty} V_x(\alpha_{nx}, \beta_{nx})$:
\begin{align*}
\limsup_{n \to \infty}  V_x(\alpha_{n x},\beta_{n x}) 
&\leq \limsup_{n \to \infty} \left\{\frac1{1-\gamma_x}{E}[\max\{0, \theta_x-c\}|~\theta_x \sim \Bta(\alpha_{n x}, ~\beta_{n x})] \right\}\\
&\leq \frac1{1-\gamma_x} \max\{0,\mu_x-c\}
\end{align*}

The second inequality holds because of the Portmanteau Theorem (\cite{Resnick2005}, page 264) since the function $0 \leq \max\{0, \theta_x-c\} \leq 1$ is bounded and continuous and the sequence of probability measures $\Bta(\alpha_{nx}, \beta_{nx})$ converges in measure to one in which $\theta_x=\mu_x$ almost surely. 


Combining the lower bound on $\liminf_{n\to \infty} V_x(\alpha_{nx}, \beta_{nx})$ and the upper bound on $ \limsup_{n\to \infty} V_x(\alpha_{nx}, \beta_{nx})$, we have 
$$ \frac 1{1-\gamma_x} \max\{0, \mu_x-c\} \le \liminf_{n\to \infty} V_x(\alpha_{nx}, \beta_{nx})
\le  \limsup_{n \to \infty} V_x(\alpha_{nx}, \beta_{nx}) \le \frac 1 {1-\gamma_x} \max\{0, \mu_x-c\}.$$
Therefore, the limit of $V_x(\alpha_{nx}, \beta_{nx})$ exists and is equal to 
$$\lim_{n \to \infty}V_x(\alpha_{nx}, \beta_{nx}) = \frac 1{1-\gamma_x}\max\{0, \mu_x-c\}.$$

\subsection*{Proof of Proposition~\ref{re:stoppingRule}}
To show the claimed characteristic of $U^*_{nx}$, it is enough to show that $U^*_{n,x} = 0$ implies $U^*_{n+1, x}=0$. To keep the notation lighter, we suppose that $n=0$, but the argument is the same for $n>0$. Let $\pi_{0u}^{(x)}$ be the policy that chooses $U_{1x}=0$ and $U_{2x}=u \in \{0, 1\}$, and then behaves optimally afterward. Suppose it is optimal to not forward at time 1, that is, $U^*_{1x}=0$ with $Q_x(\alpha_{0x}, \beta_{0x}, 1) \leq Q_x(\alpha_{0x}, \beta_{0x}, 0)$. Then, since $(\alpha_{1x}, \beta_{1x})=(\alpha_{0x}, \beta_{0x})$ under the optimal policy, either $\pi_{00}^{(x)}$ or $\pi_{01}^{(x)}$ (or both) is optimal. Suppose for contradiction $\pi_{00}^{(x)}$ is not optimal, then $\pi_{01}^{(x)}$ is strictly better than $\pi_{00}$. That is, $Q_x(\alpha_{1 x}, \beta_{1 x}, 1) > Q_x(\alpha_{1 x}, \beta_{1 x}, 0)$. But $(\alpha_{1x}, \beta_{1x})=(\alpha_{0x}, \beta_{0x})$, thus 
$$Q_x(\alpha_{0x}, \beta_{0x}, 1) = Q_x(\alpha_{1x}, \beta_{1x}, 1) > Q_x(\alpha_{1 x}, \beta_{1 x}, 0) = Q_x(\alpha_{0x}, \beta_{0x}, 0),$$
which contradicts the fact that $Q_x(\alpha_{0x}, \beta_{0x}, 1) \leq Q_x(\alpha_{0x}, \beta_{0x}, 0)$.

\subsection*{Proof of Lemma~\ref{l:V_nonDec}}
First, we show the non-decreasing property holds for a finite-horizon problem by induction. Let $V_x(\alpha, \beta ,M')$ be the value function for a problem in which we stop at time $M$ and receive terminal reward $\frac 1{1-\gamma}\max\{0, \mu_{Mx}-c\}$, and our state at time $M$ is $\alpha_{M'x} = \alpha$, $\beta_{M'x}=\beta$. 
For the base case at termination $M'=M$, we have
$$V_x(\alpha+\ell, \beta-\ell,M)=\frac 1 {1-\gamma} \max \left\{0, \frac {\alpha+\ell}{\alpha+\beta}-c \right\}$$
is indeed a non-decreasing function of $\ell$. Now suppose the property holds for some finite $M'$, with $0\le M' \le M$. That is, suppose $\ell \mapsto V_x(\alpha+\ell, \beta-\ell, M')$ is non-decreasing. Let us show that it also holds for $M'-1$. Let $$\mu(\ell)=\frac {\alpha+\ell}{\alpha+\ell+\beta-\ell}=\mu(0)+\frac \ell{\alpha+\beta}.$$ Then
\begin{align*}
V_x&(\alpha+\ell, \beta-\ell,M'-1)\\
&=\max \Big\{0, \mu(\ell)-c+\gamma_x \left[ \mu(\ell) V_x(\alpha+\ell+1, \beta-\ell, M')+ (1-\mu(\ell))V_x(\alpha+\ell, \beta-\ell+1,M')\right] \Big\}.
\end{align*}
Let $g(\ell)=V_x(\alpha+
\ell, \beta-\ell+1, M')$, where $g$ is non-decreasing by the induction hypothesis. To show that $\ell \mapsto V_x(\alpha+\ell, \beta+\ell, M')$ is non-decreasing, it is sufficient to show that
$$\ell \mapsto f(\ell)=\mu(\ell)g(\ell+1)+(1-\mu(\ell))g(\ell)$$ is non-decreasing. 
 Let $\Delta \geq 0$. Then,
\begin{align*}
f(\ell+\Delta)-f(\ell)
&=\mu(\ell+\Delta) g(\ell+\Delta+1)+(1-\mu(\ell+\Delta))g(\ell+\Delta)-\mu(\ell)g(\ell+1)-(1-\mu(\ell))g(\ell) \\
&=\left[ \mu(\ell+\Delta)-\mu(\ell) \right] g(\ell+\Delta+1) + \mu(\ell)g(\ell+\Delta+1) \\
& ~~~ +\left[ (1-\mu(\ell+\Delta))-(1-\mu(\ell))\right]g(\ell+\Delta)+(1-\mu(\ell))g(\ell+\Delta) \\
& ~~~ -\mu(\ell)g(\ell+1)-(1-\mu(\ell))g(\ell) \\
&= \left[ \mu(\ell+\Delta)-\mu(\ell)\right] [g(\ell+\Delta+1)-g(\ell+\Delta)] \\
& ~~~ +\mu(\ell) [g(\ell+\Delta+1)-g(\ell+1)]+(1-\mu(\ell))[g(\ell+\Delta)-g(\ell)]  \geq 0,
\end{align*}
since both $\mu(\ell)$ and $g(\ell)$ are non-decreasing functions of $\ell$ and $\mu(\ell) \in [0,1]$. We conclude that $V_x(\alpha+\ell, \beta-\ell, M'-1)$ is non-decreasing in $\ell$. Thus, by induction, $V_x(\alpha+\ell, \beta-\ell, M')$ is non-decreasing for all $0 \le M' \le M$. \\
Because the limit of non-decreasing functions is also non-decreasing, we have that $\ell \mapsto V_x(\alpha+\ell,\beta-\ell)=\lim_{M \rightarrow \infty} V(\alpha+\ell, \beta-\ell, M)$ is non-decreasing, as stated in Lemma \ref{l:V_nonDec}. 


\subsection*{Proof of Lemma~\ref{l:V_convex}}
Similar to the proof of Lemma~\ref{l:V_nonDec}, we show convexity of a $M$-step finite-horizon problem $V_x(\alpha, \beta, M)$ by induction, and use that the limit of convex functions is convex. Instead of $\alpha$ and $\beta$, let us rewrite the value function of the finite horizon problem in terms of $\mu=\frac{\alpha}{\alpha+\beta}$ and $m=\alpha+\beta$: 
$$f_M(\mu) \equiv V_x(\mu \cdot m, (1-\mu) \cdot m, M).$$

By induction, we show that $f_\ell(\cdot)$ is convex for all $\ell \leq M$. At the termination step, $f_M(\mu)=\frac 1 {1-\gamma_x} \max \left\{0,  \mu-c \right\}$ is convex in $\mu$. Assume that $f_{\ell+1}(\cdot)$ is convex, we show next that $f_\ell(\cdot)$ is convex. 

Let $g_\ell(\mu)=\mu-c + \gamma_x \left[ \mu \cdot f_{\ell+1} \left(\frac{\mu \ell +1}{\ell+1} \right)  + (1-\mu) \cdot f_{\ell+1} \left(\frac{\mu \ell}{\ell+1} \right) \right]$, so $f_\ell(\mu) =\max \left \{0, g_\ell(\mu)\right\}$. Since the maximum of convex functions is convex, we only need to show that $g_{\ell}(\mu)$ is convex, or equivalently, 
\begin{equation*}
g_{\ell} \left(\frac{\mu+\nu}{2} \right) \le \frac 12 (g_{\ell}(\mu)+g_{\ell}(\nu) )
\end{equation*}
 for any $0 \le \mu \le 1$ and $0 \le \nu \le 1$ with $\mu \neq \nu$, or equivalently
\begin{eqnarray*}
&&\frac{\mu+\nu}{2} -c + \gamma_x \left[\frac{\mu+\nu}{2}\cdot f_{\ell+1} \left(\frac{\frac{\mu+\nu}{2} \ell + 1}{\ell+1} \right) + \frac{2-\mu-\nu}{2}\cdot f_{\ell+1} \left(\frac{\frac{\mu+\nu}{2}\ell}{\ell+1} \right)\right] \\  
 &\le&  \frac{\mu+\nu}{2} -c +\gamma_x \left[\frac\mu2\cdot f_{\ell+1} \left(\frac{\mu \ell +1}{\ell+1}\right)  + \frac{1-\mu}{2}\cdot f_{\ell+1}\left(\frac{\mu \ell }{\ell+1}\right)
+ \frac\nu2\cdot f_{\ell+1} \left(\frac{\nu \ell +1}{\ell+1}\right) + \frac{1-\nu}2\cdot f_{\ell+1}\left(\frac{\nu \ell }{\ell+1} \right)   \right], 
\end{eqnarray*}

or equivalently, canceling some terms (first subtract $\frac{\mu+\nu}2-c$ and then divide by $\gamma_x$),
\begin{align}
\begin{split}
&\frac{\mu+\nu}{2}\cdot f_{\ell+1} \left(\frac{\frac{\mu+\nu}{2} \ell + 1}{\ell+1} \right) + \frac{2-\mu-\nu}{2}\cdot f_{\ell+1} \left(\frac{\frac{\mu+\nu}{2}\ell}{\ell+1} \right) \\
 \le & \quad  \frac\mu2\cdot f_{\ell+1} \left(\frac{\mu \ell +1}{\ell+1}\right)  + \frac{1-\mu}{2}\cdot f_{\ell+1}\left(\frac{\mu \ell }{\ell+1}\right)
+ \frac\nu2\cdot f_{\ell+1} \left(\frac{\nu \ell +1}{\ell+1}\right) + \frac{1-\nu}2\cdot f_{\ell+1}\left(\frac{\nu \ell }{\ell+1} \right).
\label{eq:convexity_l4_1_2}
\end{split}
\end{align}

Notice that by convexity
$$f_{\ell+1}\left(\frac{\frac{\mu+\nu}{2} \ell + 1}{\ell+1} \right)  \le\frac12 f_{\ell+1} \left(\frac{\mu \ell +1}{\ell+1} \right)  + \frac12 f_{\ell+1}\left(\frac{\nu \ell +1}{\ell+1} \right), \text{ and }$$
$$f_{\ell+1}\left(\frac{\frac{\mu+\nu}{2} \ell }{\ell+1} \right)  \le\frac12 f_{\ell+1} \left(\frac{\mu \ell }{\ell+1} \right)  + \frac12 f_{\ell+1}\left(\frac{\nu \ell }{\ell+1} \right),$$
then the reminder on the left--hand--side of inequality~\eqref{eq:convexity_l4_1_2}  is less than or equal to
\begin{eqnarray}
\frac{\mu+\nu}4 f_{\ell+1} \left(\frac{\mu \ell +1}{\ell+1} \right)  +\frac{\mu+\nu}4 f_{\ell+1} \left(\frac{\nu \ell +1}{\ell+1} \right) 
+ \frac{2-\mu - \nu}4 f_{\ell+1} \left(\frac{\mu \ell }{\ell+1} \right)  + \frac{2-\mu - \nu}4 f_{\ell+1} \left (\frac{\nu \ell }{\ell+1} \right).~~~~
\label{eq:convexity_l4_2}
\end{eqnarray}

Let us take the difference between the right--hand--side of \eqref{eq:convexity_l4_1_2} and \eqref{eq:convexity_l4_2} and prove that it is nonnegative. That is, let us show
\begin{eqnarray*}
&& \frac{\mu - \nu}4 f_{\ell+1}\left(\frac{\mu \ell +1}{\ell+1}\right)  -\frac{\mu-\nu}4 f_{\ell+1} \left(\frac{\nu \ell +1}{\ell+1}\right)   
- \frac{\mu - \nu}4 f_{\ell+1} \left(\frac{\mu \ell }{\ell+1}\right)  + \frac{\mu - \nu}4  f_{\ell+1}\left(\frac{\nu \ell }{\ell+1} \right)    \ge 0,
\end{eqnarray*}
or, equivalently, dividing by $\frac{\mu-\nu}4$,
\begin{eqnarray}
&&  f_{\ell+1}\left(\frac{\nu \ell }{\ell+1}\right) + f_{\ell+1}\left(\frac{\mu \ell +1}{\ell+1}\right) \ge  f_{\ell+1}\left(\frac{\mu \ell}{\ell+1}\right) +f_{\ell+1}\left(\frac{\nu \ell +1}{\ell+1}\right).
\label{eq:convexity_l4_3}
\end{eqnarray}

Without loss of generality, let us assume $\mu>\nu$.
Let $\delta=\frac{1}{(\mu-\nu)\ell+1}$ with $0 <\delta \le 1$, then we can represent,
$\frac{\mu \ell}{\ell+1} = \delta \frac{\nu \ell}{\ell+1}+ (1-\delta) \frac{\mu \ell+1}{\ell+1},$
and $\frac{ \nu \ell+1}{\ell+1} = (1-\delta) \frac{\nu \ell}{\ell+1} + \delta \frac{\mu \ell+1} {\ell+1}.$
Since $f_{\ell+1}(\cdot)$ is convex, we have
$$ \delta \cdot f_{\ell+1}\left( \frac{\nu \ell}{\ell+1} \right) + (1-\delta)\cdot f_{\ell+1}\left( \frac {\mu \ell+1} {\ell+1} \right) \ge f_{\ell+1}\left (\frac{\mu \ell}{\ell+1} \right),$$
$$ (1-\delta) \cdot f_{\ell+1} \left( \frac{ \nu \ell}{\ell+1} \right) + \delta \cdot f_{\ell+1} \left( \frac {\mu \ell+1}{\ell+1} \right) \ge f_{\ell+1}\left(\frac{\nu \ell+1}{\ell+1} \right).$$ 
Combining them, we show equality~\eqref{eq:convexity_l4_3}. 

\subsection*{Proof of Theorem~\ref{th:muStar_structure}}
\paragraph{Part 1: we want to show that $\mu^*(m) \le c$.}

As in the proof of Lemma~\ref{l:V_convex}, we write the value function in terms of $\mu=\frac{\alpha}{\alpha+\beta}$ and $m=\alpha+\beta$ instead of $\alpha, \beta$. First, $\mu V_x(\mu m+1, (1-\mu)m) + (1-\mu)V_x(\mu m, (1-\mu)m+1) \geq 0$ since both $V_x(\mu m+1, (1-\mu)m)$ and $V_x(\mu m, (1-\mu) m+1)$ are non-negative and any convex combination of two non-negative points is also non-negative. 
Pick any $\mu \in (c, 1]$, we have 
$Q_x(\mu m, (1-\mu)m, 1)=\mu-c+\gamma_x \left[\mu V_x(\mu m+1, (1-\mu)m) + (1-\mu)V_x(\mu m, (1-\mu)m+1) \right] \geq \mu-c>0, $
so $\mu \in A(m) = \Big \{\mu(\alpha, \beta): \alpha \in (0, \infty), \beta \in (0, \infty), m=\alpha+\beta \text{ and } Q_x(\alpha, \beta, 1) > 0  \Big \}.$ Therefore, $A(m) \supseteq (c, 1]$, which implies $\mu^*(m)= \inf A(m) \leq c$. \\

\paragraph{Part 2: we want to show that $\mu^*(m_0) \le \mu^*(m_0+1)$ for any $m_0>0$.}

It is sufficient to show that for each $\mu$, $V_x(\mu m_0, (1-\mu) m_0) \ge V_x(\mu (m_0+1), (1-\mu) (m_0+1))$. As defined in the proof of Lemma~\ref{l:V_convex}, we first show the statement for a $M$--step finite--horizon problem where $M-m_0\ge 1$ is an integer. Fix $\mu$, let $g_m(\mu) \equiv V_x(\mu m, (1-\mu)m, M)$ for any $0<m \le M$. We will show $g_m(\mu) \ge g_{m+1}(\mu)$ for all $m \in \{M-1, M-2, ..., m_0\}$ by backward induction. In the base case, $m=M-1$, we have
\begin{eqnarray*}
&g_{M-1}(\mu) &= \max \left \{ 0, \mu - c + \gamma_x \left[ \mu \cdot g_{M} \left( \frac{\mu (M-1)+1}{M} \right)+ (1-\mu) \cdot g_{M} \left(\frac{\mu (M-1)}{M} \right)\right] \right\} \\
&&\ge  \max \left \{ 0, \mu - c + \gamma_x \cdot g_M(\mu) \right\} =\max \left\{ 0, \mu-c+\gamma_x [\max\{0, (\mu-c)/(1-\gamma_x)\}] \right\},
\end{eqnarray*}
where the inequality holds by Jensen's inequality because $g_M(\mu)=\max\{0, \mu-c\}/(1-\gamma_x)$ is convex in $\mu$. If $\mu-c>0$, then $g_{M-1}(\mu)\ge \frac {\mu-c}{1-\gamma_x}=g_M(\mu)$. Otherwise, if $\mu-c \le 0$, then $g_{M-1}(\mu) \ge 0 = g_M(\mu)$. This shows the base case $g_{M-1}(\mu)\ge g_M(\mu)$.

Let $m \in \{M-2, M-3, ..., m_0\}$ and assume $g_{m}(\mu) \ge g_{m+1}(\mu)$. Next we show that $g_{m-1}(\mu) \ge g_{m}(\mu)$. There are two cases to consider. First, if $g_m(\mu)=0$, we have that $g_{m-1}(\mu) \ge g_m(\mu)= 0$. Second, if

\begin{center} 
$g_m(\mu)= \mu-c +\gamma_x \left[\mu \cdot g_{m+1} \left( \frac{ \mu m+1}{m+1} \right)+ (1-\mu) \cdot g_{m+1} \left( \frac{ \mu m}{m+1} \right) \right],$
\end{center}

\noindent then
\begin{eqnarray*}
&&\frac{g_{m-1}(\mu) - g_{m}(\mu)}{\gamma_x} \\
 &\ge& \mu \cdot g_m \left(  \frac{\mu (m-1) +1}{m} \right) + (1-\mu) \cdot g_m  \left( \frac{\mu (m-1)}{m} \right) -\mu \cdot g_{m+1} \left( \frac{ \mu m+1}{m+1} \right)- (1-\mu) \cdot g_{m+1} \left( \frac{ \mu m}{m+1} \right) \\
& \ge & \mu \cdot g_m \left(\frac{\mu (m-1) +1}{m} \right) + (1-\mu) \cdot g_m  \left( \frac{\mu (m-1)}{m} \right) 
 -\mu \cdot g_{m} \left( \frac{ \mu m+1}{m+1} \right)- (1-\mu) \cdot g_{m} \left( \frac{ \mu m}{m+1} \right),
\end{eqnarray*} 
in which the first inequality is because $g_{m-1}$ is greater than or equal to its value at forwarding, and the second inequality is due to the induction hypothesis that $g_{m+1}(\mu) \le g_{m}(\mu)$. Notice that 
$\frac{\mu(m-1)}{m} \le \frac{\mu m}{m+1} \le \mu \le \frac{\mu m+1}{m+1} \le \frac{\mu(m-1)+1}{m}.$ 
Define a function $\hat g_m(x)$ by replacing $g_m(x)$ with a linear interpolation between $g_m \left(\frac{\mu m}{m+1} \right)$ and $g_m \left( \frac{\mu m+1}{m+1} \right)$ over the interval between these two points, 
\begin{equation} 
\hat g_m (x)= \left\{ \begin{array} {ll} g_m(x), & \mbox{ ~~if } x\ge \frac{\mu m+1}{m+1}  \mbox{ or } x \le \frac{\mu m}{m+1}, \\  \medskip\\
\frac{g_m \left(\frac{\mu m+1}{m+1} \right)-g_m \left( \frac{\mu m}{m+1} \right)}{ 1/(m+1)}   \left (x-\frac{\mu m}{m+1} \right)  + g_m \left(\frac{\mu m}{m+1} \right),& \mbox{ ~~if } \frac{\mu m}{m+1}\le  x \le \frac{\mu m+1}{m+1}. 
\end{array}\right.
\end{equation}

$g_m$ is convex by Lemma~\ref{l:V_convex}, and so $\hat g_m$ is also convex. Moreover, by Jensen's inequality, $\hat g_m \ge g_m$. Thus, using this inequality in the first line, and then Jensen's inequality on $\hat g_m$ in the second line, 
\begin{eqnarray*}
(1-\mu)\cdot g_m \left(\frac{\mu (m -1)}{m} \right) &+&\mu \cdot g_m \left( \frac{\mu(m -1)+1}{m} \right)
\ge (1-\mu) \cdot \hat g_m \left(\frac{\mu (m -1)}{m} \right) +\mu \cdot \hat g_m \left( \frac{\mu(m -1)+1}{m} \right)\\
&\ge & \hat g_m \left((1-\mu)\frac{\mu (m -1)}{m}  +\mu \frac{\mu (m -1)+1}{m}   \right)
\\
&=&\hat g_m(\mu) = \hat g_m \left((1-\mu)\frac{\mu m}{m+1}+ \mu \frac{\mu m+1}{m+1} \right )\\
&=& \frac{g_m \left(\frac{\mu m+1}{m+1} \right)-g_m \left( \frac{\mu m}{m+1} \right)}{ 1/(m+1)}   \left(\mu-\frac{\mu m}{m+1} \right)  + g_m \left(\frac{\mu m}{m+1} \right)
\\
&=& (1-\mu) \cdot g_m \left(\frac{\mu m}{m+1}\right)+ \mu \cdot g_m \left(\frac{\mu m+1}{m+1}\right).
\end{eqnarray*}
Thus, $g_{m-1}(\mu) \ge g_m(\mu)$. This completes the induction step. 

We have shown that $V_x(\mu m_0, (1-\mu)m_0, M) = g_m(\mu)\ge V_x(\mu (m_0+1), (1-\mu)(m_0+1), M)$ for each $\mu$ in a $M$--step finite--horizon problem. Taking the limit as $M$ goes to infinity, we conclude that $V_x(\mu m_0, (1-\mu)m_0) \ge V_x(\mu (m_0+1), (1-\mu)(m_0+1))$. Equivalently, we conclude that $\mu^*(m) \le \mu^*(m+1)$. \\

\paragraph{Part 3: we want to show that $\lim_{m \rightarrow \infty} \mu^*(m)=c$.}

It is sufficient to show that, for each $\mu < c$, there exists $N>0$ large enough that $V_x(\mu m, (1-\mu)m)=0$ for all $m \ge N$. Let $\alpha = \mu m$, and $\beta=(1-\mu)m$. By definition,
$V_x(\alpha, \beta) = \max\{ 0, \mu-c+\gamma_x [ \mu V_x(\alpha+1, \beta) + (1-\mu) V_x(\alpha, \beta+1)]\},$
so if $\mu-c+\gamma_x [ \mu V_x(\alpha+1, \beta) + (1-\mu) V_x(\alpha, \beta+1)] \le 0$, then $V_x(\alpha, \beta)=0$. Using the upper bound from Proposition~\ref{prop:upperBound}, $V_x(\alpha, \beta) \le V^U_x(\alpha, \beta)= \frac 1{1-\gamma_x} E [ \max \{0, \theta_x-c\} | \theta_x \sim \Bta(\alpha, \beta)]$. So,
\begin{align*}
\mu &V_x(\alpha+1, \beta) + (1-\mu) V_x(\alpha, \beta+1)\le \mu V^U_x(\alpha+1, \beta) + (1-\mu) V^U_x(\alpha, \beta+1)\\
&=\mu E[\max \{0, \theta_x-c\}|\theta_x \sim \Bta(\alpha+1, \beta)] + (1-\mu)E[\max \{0, \theta_x-c\}|\theta_x \sim \Bta(\alpha, \beta+1)]\\
&=E[E[\max \{0, \theta_x-c\}|\theta_x \sim \Bta(\alpha, \beta), Y]|\theta_x \sim \Bta(\alpha, \beta)]\\
&=E[\max\{0, \theta_x-c\}|\theta_x \sim \Bta(\alpha, \beta)]
=V^U_x(\alpha, \beta),
\end{align*}
where the first and fourth equalities are justified through the definition of $V^U_x(\alpha, \beta)$, the second equality rewrites each term in terms of a conditional expectation given $Y|\theta_x \sim \text{Bernoulli}(\theta_x)$, and the third equality is due to the tower property of conditional expectation. Then, recalling that $\alpha$, $\beta$ implicitly depend on $m$, 
\begin{align*}
\lim_{m \rightarrow \infty} V^U_x(\alpha, \beta) 
&=\lim_{m \rightarrow \infty} E [ \max\{0, \theta_x-c\} | \theta_x \sim \Bta(\alpha, \beta) ]=\max \{0, \mu-c\} = 0.
\end{align*}
So, $\limsup_{m \rightarrow \infty} \mu-c+\gamma_x [\mu V_x(\alpha+1, \beta)+(1-\mu)V_x(\alpha, \beta+1)] \le \mu-c < 0$. Thus, there exists $N>0$ such that $\mu-c+\gamma_x(\mu V_x(\alpha+1, \beta) + (1-\mu) V_x(\alpha, \beta+1)] <0$ for all $m \ge N$. This implies that $\mu^*(m) \ge \mu$ $\forall m > N$. Then,
$\liminf_{n \rightarrow \infty} \mu^*(m) = \lim_{N' \rightarrow \infty} \inf \{\mu^*(m): m \ge N'\} \ge \mu.$
Since this is true for all $\mu<c$, we have $\liminf_{m \rightarrow \infty} \mu^*(m) \ge c.$ Combining with part 1, which showed that $\limsup_{m \rightarrow \infty} \mu^*(m) \le c$, we have that 
$$c \le \liminf_{m \rightarrow \infty} \mu^*(m) \le \limsup_{m \rightarrow \infty} \mu^*(m) \le c.$$
Therefore, the limit exists and $\lim_{m \rightarrow \infty} \mu^*(m)=c$. 

\subsection*{Proof of Lemma~\ref{lemma:bounds}}
At each time step, the stepwise reward function $U_{nx}(Y_{nx}-c)$ is bounded above by $1$.  Then we have
\begin{eqnarray*}
V_x(\alpha, \beta) 
&= \sup_{\pi} E^{\pi} \left [ \sum_{n=1}^\infty \gamma_x^{n-1}(Y_{nx}-c)U_{nx} ~|~  \theta_x \sim \Bta(\alpha, \beta) \right ] \leq \sup_{\pi} E^{\pi} \left [ \sum_{n=0}^\infty \gamma_x \right] = \frac {1}{1-\gamma_x}.
\end{eqnarray*} 
For the other side, any policy $\pi'$ provides a lower bound on the value of an optimal policy.  Thus,
\begin{align*}
V_x(\alpha, \beta) 
&=\sup_{\pi} E^{\pi} \left [ \sum_{n=1}^\infty \gamma_x^{n-1}(Y_{nx}-c)U_{nx} ~|~  \theta_x \sim \Bta(\alpha, \beta) \right ] 
\ge E^{\pi'} \left [ \sum_{n=1}^\infty \gamma_x^{n-1}(Y_{nx}-c)U_{nx} ~|~  \theta_x \sim \Bta(\alpha, \beta) \right ].
\end{align*}
Take $\pi'$ to be the (deterministic) policy that chooses, for each $n$, $U_{nx}=1$ if $\mu(\alpha, \beta)>c$, and $U_{nx}=0$ if $\mu(\alpha, \beta)\leq c$. 
If $\mu(\alpha,\beta)-c>0$, its value is
\begin{equation*}
E^{\pi'} \left [ \sum_{n=1}^\infty \gamma_x^{n-1}(Y_{nx}-c) ~|~  \theta_x \sim \Bta(\alpha, \beta) \right ]
= \sum_{n=1}^\infty \gamma_x^{n-1}(\mu(\alpha,\beta) -c)
= \frac1{1-\gamma_x} (\mu(\alpha,\beta) -c).
\end{equation*}
If $\mu(\alpha,\beta)-c\le 0$, its value is $0$.
Thus, putting both cases together, we see that its value is 
\begin{equation*}
E^{\pi'} \left [ \sum_{n=1}^\infty \gamma_x^{n-1}(Y_{nx}-c) ~|~  \theta_x \sim \Bta(\alpha, \beta) \right ]
=\frac1{1-\gamma_x} \max\{ \mu(\alpha,\beta) -c, 0 \}.
\end{equation*}
Combining this expression with the previously stated lower bound on $V_x(\alpha,\beta)$ shows that 
$V_x(\alpha,\beta) \ge \frac1{1-\gamma_x} \max\{ \mu(\alpha,\beta) -c, 0 \}$.

\subsection*{Proof of Lemma~\ref{l:approximationBounds}}
Consider the base case at termination $M$.  By Lemma \ref{lemma:bounds}, we have 
$$\frac1{1-\gamma_x} \max \left\{0, \frac{\alpha}{\alpha+\beta}-c \right\} = V^L_x(\alpha, \beta; M) \leq V_x(\alpha, \beta) \leq V_x^U(\alpha, \beta;M) = \frac 1{1-\gamma_x}.$$
Assume that the statement holds at some $0\leq \ell+1\leq M$, we show that it also holds at $\ell$: 
\begin{align*}
V_x^U(\alpha, \beta; \ell) 
&= \max \left \{0, \frac{\alpha}{\alpha+\beta}-c+\gamma_x\left[ \frac{\alpha}{\alpha+\beta} V^U_x(\alpha+1, \beta; M) + \frac{\alpha}{\alpha+\beta}V_x^U(\alpha, \beta+1;M) \right] \right\} \\
& \geq \max \left\{0, \frac{\alpha}{\alpha+\beta}-c+\gamma_x\left[ \frac {\alpha}{\alpha+\beta} V_x(\alpha+1, \beta; M) +\frac{\alpha}{\alpha+\beta}V_x(\alpha, \beta+1;M) \right] \right\} =V_x(\alpha, \beta)
\end{align*}
Similarly, 
\begin{align*}
V_x^L(\alpha, \beta; M) 
&= \max \left \{0, \frac{\alpha}{\alpha+\beta}-c+\gamma_x\left[ \frac{\alpha}{\alpha+\beta} V^L_x(\alpha+1, \beta; M) + \frac{\alpha}{\alpha+\beta}V_x^L(\alpha, \beta+1;M) \right] \right\} \\
& \leq \max \left\{0, \frac{\alpha}{\alpha+\beta}-c+\gamma_x\left[ \frac {\alpha}{\alpha+\beta} V_x(\alpha+1, \beta; M) +\frac{\alpha}{\alpha+\beta}V_x(\alpha, \beta+1;M) \right] \right\} =V_x(\alpha, \beta)
\end{align*}
Combining the upper bound and lower bound, we show that $V^L_x(\alpha, \beta;M) \leq  V_x(\alpha, \beta) \leq V^U_x(\alpha, \beta;M)$ for each $M\geq 0$ and for each $(\alpha, \beta)$ pair.

\subsection*{Proof of Lemma~\ref{l:gap_convergence}}
For any $M \geq 0$, let us first consider the base case at termination $M$. Then, by definition, for any $0\leq i \leq M$ with $\alpha=\alpha_{0x}+i$ and $\beta=\beta_{0x}+M-i$,  we have $V^U_x(\alpha, \beta; M)=\frac 1 {1-\gamma_x}$ and $V^L_x(\alpha, \beta; M)=\frac 1{1-\gamma_x} \max \left \{0, \frac{\alpha}{\alpha+\beta}-c \right\}\geq 0$, thus, $V^U_x(\alpha, \beta; M)-V^L_x(\alpha, \beta; M) \leq \frac 1 {1-\gamma_x}$. 

Assume that the statement holds for some $0\le \ell+1 \leq M$. We show in the following that the inequality also holds for $\ell$ by induction. There are two cases to consider at step $\ell$: 

Case 1: Suppose $V^U_x(\alpha, \beta; M)=0$. Since $V^L_x(\alpha, \beta; M) \geq 0$, the induction step at $\ell$ is obviously true. 

Case 2: Suppose $V^U_x(\alpha, \beta; M)>0$. Then, 
$$V^U_x(\alpha, \beta;M) = \frac {\alpha}{\alpha+\beta}-c+\gamma_x \left[\frac {\alpha}{\alpha+\beta} V^U_x(\alpha+1, \beta; M)+\frac{\beta}{\alpha+\beta} V^U_x(\alpha, \beta+1;M) \right], \text{ and }$$
$$ V^L_x(\alpha, \beta;M) \geq \frac {\alpha}{\alpha+\beta}-c+\gamma_x \left[\frac{\alpha}{\alpha+\beta} V^L_x(\alpha+1, \beta;M)+\frac{\beta}{\alpha+\beta}V^L_x(\alpha, \beta+1;M) \right].$$
Combining them leads to 
\begin{align*}
V^U_x(\alpha, \beta;M)-V^L_x(\alpha, \beta;M)
&\leq \gamma_x \left[\frac{\alpha}{\alpha+\beta} \left[V^U_x(\alpha+1, \beta;M)-V^L_x(\alpha+1, \beta;M)\right] \right. \\
& ~~~~~~~~+ \left. \frac{\beta}{\alpha+\beta} \left[ V^U_x(\alpha, \beta+1;M)-V^L_x(\alpha, \beta+1;M) \right] \right] \\
& \leq \gamma_x \left[ \left(\frac{\alpha}{\alpha+\beta}  + \frac{\beta}{\alpha+\beta} \right) \frac{\gamma_x^{M-(\ell+1)}}{1-\gamma_x} \right] = \frac{\gamma_x^{M-\ell}}{1-\gamma_x},
\end{align*}
while the second inequality is due to the inductive assumption that the statement holds for $\ell+1$, that is, $V^U_x(\alpha+1, \beta;M)-V^L_x(\alpha+1, \beta;M) \leq \frac{\gamma_x^{M-(\ell+1)}}{1-\gamma_x}$ and $V^U_x(\alpha, \beta+1;M)-V^L_x(\alpha, \beta+1;M) \leq \frac{\gamma_x^{M-(\ell+1)}}{1-\gamma_x}$. Therefore, we conclude that it holds for any $0 \leq \ell \leq M$. Moreover, as $M \rightarrow \infty$, the difference between $V^U_x(\alpha, \beta; M)$ and $V^L_x(\alpha, \beta; M)$ converges to $0$ for any $0< \gamma_x < 1$. 


\bibliographystyle{ormsv080}
\bibliography{xzhao}
\end{document}